\def\draft{n}
\documentclass{amsart}
\usepackage{fullpage,amssymb,epic,eepic,epsfig}

\theoremstyle{plain}

\newtheorem{theorem}{Theorem}
\newtheorem{observation}{Observation}
\newtheorem{proposition}{Proposition}[section]
\newtheorem{lemma}[proposition]{Lemma}
\newtheorem{corollary}[proposition]{Corollary}

\theoremstyle{definition}
\newtheorem{definition}[proposition]{Definition}

\theoremstyle{remark}
\newtheorem{example}[proposition]{Example}

\newtheorem{remark}[proposition]{Remark}

\def\printname#1{
        \if\draft y
                \smash{\makebox[0pt]{\hspace{-0.5in}
                        \raisebox{8pt}{\tt\tiny #1}}}
        \fi
}

\newcommand{\psdraw}[2]
         {\begin{array}{c} \hspace{-1.3mm}
        \raisebox{-4pt}{\epsfig{figure=draws/#1.eps,width=#2}}
        \hspace{-1.9mm}\end{array}}

\newlength{\standardunitlength}
\catcode`\@=11
\long\def\@makecaption#1#2{%
     \vskip 10pt

\setbox\@tempboxa\hbox{
       \small\sf{\bfcaptionfont #1. }\ignorespaces #2}%
     \ifdim \wd\@tempboxa >\captionwidth {%
         \rightskip=\@captionmargin\leftskip=\@captionmargin
         \unhbox\@tempboxa\par}%
       \else
         \hbox to\hsize{\hfil\box\@tempboxa\hfil}%
     \fi}
\font\bfcaptionfont=cmssbx10 scaled \magstephalf
\newdimen\@captionmargin\@captionmargin=2\parindent
\newdimen\captionwidth\captionwidth=\hsize
\catcode`\@=12

\newcommand{\tr}{\operatorname{tr}}

\def\lbl#1{\label{#1}\printname{#1}}

\def\BN{\mathbb N}
\def\BZ{\mathbb Z}

\def\A{\mathcal A}

\def\C{\mathcal C}

\def\D{\Delta}
\def\K{\mathcal K}

\def\F{\mathcal F}

\def\P{\mathcal P}
\def\R{\mathcal R}

\def\I{\mathcal I}
\def\CI{\mathcal C}
\def\AC{\mathcal{AC}}

\def\Ga{\Gamma}

\def\ga{\gamma}

\def\w{\omega}

\def\e{\epsilon}
\def\Ga{\Gamma}
\def\d{\delta}
\def\b{\beta}

\def\sub{\subseteq}
\def\inv{\mathrm{inv}}

\def\ti{\widetilde}

\def\bt{\bar{t}}
\def\rot{\mathrm{rot}}
\def\longto{\longrightarrow}
\def\Sev{{\mathcal S}}

\def\sign{\mathrm{sign}}
\def\exc{\mathrm{exc}}

\def\mult{\mathrm{mult}}

\def\smooth{\mathrm{smooth}}
\def\bq{\bar{q}}
\def\err{\mathrm{err}}
\def\dec{\mathrm{dec}}
\def\vert{\mathrm{vert}}
\def\circ{\mathrm{circ}}
\def\longto{\longrightarrow}
\def\unknot{\mathrm{unknot}}
\def\defect{\mathrm{def}}

\def\Edges{\mathrm{Edges}}

\def\Ferm{\mathrm{Ferm}}

\begin{document}


\title[A non-commutative formula for the colored Jones function]
{A non-commutative formula for the colored Jones function}
\author{Stavros Garoufalidis}
\address{School of Mathematics \\
         Georgia Institute of Technology \\
         Atlanta, GA 30332-0160, USA}
\email{stavros@math.gatech.edu}
\author{Martin Loebl}
\address{KAM MFF UK \\
Institute of Theoretical Computer Science (ITI)\\
Charles University \\
Malostranske n. 25 \\
118 00 Praha 1 \\
Czech Republic.}
\email{loebl@kam.mff.cuni.cz}

\thanks{S.G. was supported in part by National Science Foundation. 
M.L. gratefully acknowledges the support of ICM-P01-05. This work was done 
partly while M.L. was visiting DIM, U. Chile.
       \newline
1991 {\em Mathematics Classification.} Primary 57N10. Secondary 57M25.
}

\date{
This edition: March 7, 2005 \hspace{0.5cm} First edition: November 22, 2004.}

\begin{abstract}
The colored Jones function of a knot is a sequence of Laurent polynomials
that encodes the Jones polynomial of a knot and its parallels. It has been 
understood 
in terms of representations of quantum groups and Witten gave an intrinsic 
quantum field theory interpretation of the colored Jones function as 
the expectation 
value of Wilson loops of a 3-dimensional gauge theory, the Chern-Simons 
theory. 
We present the colored Jones function as an evaluation of the inverse of a 
non-commutative 
fermionic partition function. This result is in the form familiar in 
quantum field theory, namely the inverse of a generalized determinant.
 Our formula also reveals a direct relation between the Alexander
polynomial and the colored Jones function of a knot and immediately implies
 the extensively studied Melvin-Morton-Rozansky conjecture,
first proved by Bar-Natan and the first author about ten years ago.
\end{abstract}

\maketitle

\tableofcontents


\section{Introduction}
\lbl{sec.intro}

\subsection{The Jones polynomial of a knot}
\lbl{sub.history}
In 1985, V. Jones discovered a celebrated invariant of knots, the Jones
polynomial, \cite{J1}. 
Jones's original formulation of the Jones polynomial was given in terms of
{\em representations of braid groups} and {\em Hecke algebras}, \cite{J1}. 
It soon
became apparent that the Jones polynomial can be defined as a {\em state sum} 
of a {\em statistical mechanics} model that uses as input a planar projection
of a knot, \cite{J2,Tu}.
As soon as the Jones polynomial was discovered, it was compared 
with the better-understood Alexander polynomial of a knot. The latter can be
defined using classical algebraic topology (such as the homology of 
the infinite cyclic cover of the knot complement), and its skein theory
can be understood purely topologically.
On the other hand, the Jones polynomial appears to be difficult to understand
topologically, and there is a good reason for this, as was explained by 
Witten, \cite{Wi}.
Namely, the Jones polynomial can be thought of as the expectation 
value of Wilson loops of a 3-dimensional gauge theory, the Chern-Simons theory;
 in general, this is hard to understand. Witten's approach leads to a number
of conjectures that relate limits of the Jones polynomial to geometric 
invariants of a knot, such as representations of the fundamental group
of its complement into compact Lie groups. 
A recent approach to the Jones polynomial 
in terms of $D$-{\em modules} and {\em holonomic functions} seems to relate 
well 
to the hyperbolic geometry of knot complements, \cite{GLe, Ga2} and
yet another approach to the Jones polynomial is via the {\em Kauffman bracket
skein theory}, \cite{KL}.

The goal of our paper is to present the colored Jones function
as {\em an evaluation of the inverse of a non-commutative fermionic partition 
function}.
 This result is in the form very familiar in quantum field theory, 
namely the inverse of
a generalized determinant.
Hence there should be a quantum field theoretic derivation of it, which may 
teach us new things about how to compute path integrals in topological 
quantum field theory.

About 10 years ago, Melvin-Morton and Rozansky independently conjectured
a relation among the limiting behavior of the colored Jones function 
of a knot and its Alexander polynomial (see Corollary \ref{cor.MMR}), 
\cite{MM,Ro1,Ro2}. D. Bar-Natan and 
the first author reduced the conjecture about knot invariants to a statement 
about their {\em combinatorial weight systems}, and then proved it for all 
weight systems that come from {\em semisimple Lie algebras} using 
combinatorial 
Lie algebraic methods, \cite{BG}. A combinatorial description of the 
corresponding
weight systems was obtained in \cite{GL}. Over the years, the MMR Conjecture 
has received 
attention by many researchers who gave alternative proofs, 
\cite{Ch,KSA,KM,Ro3,V}.

A comparison of Theorem \ref{thm.al} and Theorem \ref{thm.nonc} reveals a 
direct
relation between the Alexander polynomial and the colored 
Jones function. This should help us better understand the topological features
of the colored Jones function.

We will introduce an auxiliary weighted directed graph, the arc-graph, 
that encodes transitions of walks along a planar projection of a knot.
Our results are obtained by studying the non-negative integer flows on this 
arc-graph
and applying the recently discovered  q-MacMahon Master Theorem of \cite{GLZ}.

\subsection{Statement of the main result}
\lbl{sub.sta}

\begin{definition}
\lbl{def.com}
We consider $5r$ indeterminates $r^-_i,r^+_i,u^-_i,u^+_i,z_i, 1\leq i\leq r$.
 Let $A$ be a $r$ by $r$ matrix where each indeterminate appears at most once
in an entry, and each entry is an indeterminate times a power of $q$.
 We assume $q$ is an indeterminate which commutes with all other 
indeterminates. 
Moreover we assume that each column contains at most one $u$ indeterminate, 
in its first
or last entry different from $z$. Let $L(A)$ be the set of those columns of $A$
where $u$ appears in the last non-$z$-entry. 

We define a {\em noncommutative algebra} $\A(A)$ generated by the 
indeterminates which appear in $A$, modulo 
the commutation relations specified below. Consider any $2$ by $2$ minor
 of $A$ consisting of rows $i$ and $i'$, and columns $j$ and $j'$ 
(where $1\leq i <i'\leq r$, and $1\leq j <j'\leq r$), writing $a=a_{ij}$, 
$b=a_{ij'}$, $c=a_{i'j}$, $d=a_{i'j'}$, we have the following commutation 
relations
(we will use the symbol $=_q$ to denote 'equality up to a power of $q$'): 
\begin{enumerate}
\item
The commutation in each row:
 $ba= q^{-2}ab$ if $b=_q u^-$ or $a=_q u^-$ and $ba= ab$ otherwise. 
The same rule is adapted for $cd$ commutation.
\item
The $bc$ commutation: $bc= q^{-1+s}cb$ if
 
$c=_q u^s, b=_q r^{-+}$ or 

$b=_q u^s, c=_q r^{-+}$ or

$c=_q u^s, b=_q u^{-+}, d=_q r^{-+}, a= z$ or

$b=_q u^s, c=_q u^{-+}, a=_q r^{-+}, d= z$.

$bc= q^{-1+s+s'}cb$ if $b=_q u^s, c=_q u^{s'}, a=_q r^{-+}, d=_q r^{-+}$, and

$bc = q^{-1}cb$ otherwise.
\item
Finally we require that $A$ is right-quantum (see \cite{GLZ}), i.e.,
$$ 
ca=qac, \qquad  db=qbd, \qquad 
ad=da +q^{-1}cb -qbc.
$$
\end{enumerate}
\end{definition}

Note that the commutation relations are such that
each monomial in $\A(A)$ can be brought into a q-combination of {\em canonical
monomials} $\prod_{i=1}^r a_{1i}^{m_{1i}}\dots a_{ri}^{m_{ri}}$.

\begin{definition}
\lbl{def.eval}
We define {\em n-evaluation} of a canonical monomial 
$\prod_{i=1}^r a_{1i}^{m_{1i}}\dots a_{ri}^{m_{ri}}$
to be zero if there is $ij$ with $m_{ij}> 0$ and $a_{ij}= z_k$ and otherwise
$$
\tr_n\prod_{i=1}^r a_{1i}^{m_{1i}}\dots a_{ri}^{m_{ri}}=
\prod_{i\notin L(A)}\tr_n a_{1i}^{m_{1i}}\dots a_{ri}^{m_{ri}}
\prod_{i\in L(A)}\tr_n a_{ri}^{m_{ri}}\dots a_{1i}^{m_{1i}},
$$
and
$$
\tr_n (u^{s_0})^{p_0}(r_{i_1}^{s_1})^{p_1}\dots (r_{i_m}^{s_m})^{p_m}=
q^{-s_0p_0n}\prod_{i=1}^m\prod_{j=0}^{p_i-1}
(1-t^{-s_i(n-j-p_0-\dots p_{i-1})}).
$$
\end{definition}

We consider a generic planar projection $\K$ of an oriented zero framed knot 
with $r+1$ crossings and with no kinks, together with a special arc decorated
 with $\star$. Let $K$ denote the corresponding {\em long knot} obtained by 
breaking the special arc. We will order the arcs of $\K$ so that they appear
in increasing order as we walk in the direction of the knot, such that
the special arc is last. We will also order the crossings of $\K$
such that arc $a_i$ ends at the $i$th crossing, for $i=1,\dots,r+1$.

Note that $\K$ can be uniquely reconstructed from $K$, so that any invariant
of knots gives rise to a corresponding invariant of long knots.
We consider {\em transitions} of $K$: when we walk along $K$, 
we either go under a crossing (blue transition), or jump up at a crossing
(red transition). Each transition from arc $a_i$ to arc $a_j$
is naturally equipped with a non-negative integer $\rot(a_i,a_j)$ which
can be seen from $K$ (see Definition \ref{def.jjj}).

We define the $r$ by $r$ transition matrix $B_K= (b_{ij})$ as follows.

\begin{definition}
\lbl{def.transm}
$$
b_{ij}=
\begin{cases}
q^{-\rot(ij)}u^{\sign(i)}_i & \text{if} \quad j=i+1 \\
q^{-\rot(ij)}r^{\sign(i)}_i & \text{if} \quad a_ia_j 
 \text{is a red transition} \\
z_i & \quad \text{otherwise}
\end{cases}
$$
\end{definition}
 
The next well-known theorem (see e.g. \cite{BG}) identifies the Alexander 
polynomial $\D(\K)$ 
of a knot diagram $\K$ with the determinant of $B_K$.

\begin{theorem}
\lbl{thm.al}
For every knot diagram $\K$ we have:
$$
\D(\K, t)=_t \det(I-B_K)|_{q=1, z_i=0,
u^{\sign(i)}_i=t^{-\sign(i)}, r^{\sign(i)}_i=(1-t^{-\sign(i)})}.
$$
\end{theorem}

\begin{definition}
\lbl{def.ferm}
The quantum determinant of an $r$ by $r$ matrix $A= (a_{ij})$, introduced 
in \cite{FRT},
 may be defined by
$$
\det_q(A)=\sum_{\pi\in S_r} (-q)^{-\inv(\pi)}a_{\pi(1)1}a_{\pi(2)2}\dots 
a_{\pi(r)r},
$$
where $\inv(\pi)$ equals the number of pairs $1\leq i<j \leq r$ for which 
$\pi(i)>\pi(j)$.
Moreover we let
$$
\Ferm(A)= \sum_{J\subset \{1,\dots, r\}} (-1)^{|J|}\det_q(A_J)
$$
where $A_J$ is the $J$ by $J$ submatrix of $A$.
\end{definition}

If $q=1$ then $\Ferm(A)= \det(I-A)$. Recall the MacMahon Master Theorem 
(\cite{MMM}),
known also as the {\em boson-fermion} correspondence 
$$
\frac{1}{\det(I-A)}= \sum_{n=0}^{\infty} \tr S^n(A),
$$
where $S^n(A)$ is the n-th symmetric power of $A$.

The main result of this paper is as folows:

\begin{theorem}
\lbl{thm.nonc}
For every knot diagram $\K$ we can construct a matrix $B'_K$ from $B_K$
by a permutation of rows and columns so that 
$$
J_n(\K,q) = q^{\d(K,n)}1/\Ferm(B'_K),
$$
n-evaluated; $\d(K,n)$ is an integer that can be computed easily from $K$ 
(see Definition
\ref{def.jj}). 
\end{theorem}

As an immediate consequence we obtain the seminal Melvin-Morton-Rozansky 
Conjecture (MMR in short),
whose proof was first given by \cite{BG}.

\begin{corollary}
\lbl{cor.MMR}
$$
\lim_{n \to \infty} J_n(\K,q^{1/n})=_t \frac{1}{\D(\K,q)}.
$$
\end{corollary}

\begin{remark}
The computational complexity of the Jones polynomial and its approximation
is studied extensively and as far as we know, this cannot be said about 
non-commutative formulas. Hence it may be enlightening to study our formula
from a computation point of view.
\end{remark}

\subsection{Acknowledgement} The authors wish to thank A. Kricker, 
M. Klazar, TTQ. Le and M. Staudacher for helpful discussions.
TTQ. Le has informed us that he also has a non-commutative expression
for the colored Jones function for braids, obtained by a different method. 

\section{The zeta function of a graph and the quantum MacMahon Master Theorem}
\lbl{sec.zeta}

One of main ingredients in our result is combinatorics of non-negative
integer flows on digraphs. They appear in an expression of the zeta function.

Let us recall what is the zeta function of a digraph. We will consider
digraphs (that is, graphs with oriented edges) with weights on their edges.

Let $G=(V,E)$ be a digraph with vertex set $V$ and directed edges $E \subset
V \times V$, and let $B=(\b_e)_{e\in E}$ be a weight matrix for the edges of 
$G$.
For edge $e$ we denote by $s(e), t(e)$ the starting and terminal vertex of $e$.
Bass-Ihara-Selberg defined
a zeta function of a graph in analogy with number theory and dynamical
systems, where the analogue of a prime number is a nonperiodic cycle.
Let us define the latter.

A {\em pointed walk} on a digraph is a 
sequence $(e_1, \dots, e_k)$ of edges such that the end of one coincides with 
the beginning of the next; we say that it is pointed at the beginning of $e_1$,
which is also called a {\em base point}. A {\em pointed closed walk} 
is a path whose beginning and end vertex coincide. 
Two pointed closed walks are {\em equivalent} if they differ on the choice
of base point only. By a {\em cycle} we will mean an equivalence class 
of pointed closed walks. A  cycle $c$ is {\em periodic} if $c=d^n$ for 
some closed walk $d$ and some integer $n>1$. Otherwise, it is called
{\em nonperiodic}. Let $\P(G)$ denote the set of {\em nonperiodic cycles}
of a digraph $G$. Using the weight function, we may define the weight $\b(c)$
of a cycle $c$ by $\b(c)= \prod_{e\in c}\b(e)$.

With the above preliminaries, Bass-Ihara-Selberg \cite{B} define

\begin{definition}
\lbl{def.zeta}
The {\em zeta function} $\zeta(G,B)$ of a weighted digraph $(G,B)$ is defined
by:
$$
\zeta(G,B)=\prod_{c \in \P(G)} \frac{1}{1-\b(c)}.
$$
\end{definition}

It follows by definition that
$$
\zeta(G,B)= \sum_{c \, \text{multisubsets of} \, \P(G)} \b(c)
$$

The actual definition of Bass-Ihara-Selberg uses more special
weights for the edges (each edge is given the same weight), and is used
to digraphs which are doubles (in the sense of replacing an unoriented edge
by a pair of oppositely oriented edges) of undirected graphs.
 
Foata-Zeilberger proved that the zeta function is a
rational function, and in fact given by the inverse of a determinant. Moreover,
the zeta function is given by a sum over flows. 

\begin{definition}
\lbl{def.flow}  
A {\em flow} $f$ on a digraph $G$ is a function $f:
\mathrm{Edges}(G)\longto \BN$
of the edges of $G$ that satisfies the (Kirkhoff) {\em conservation law}
$$
\sum_{e \, \text{begins at} \, v} f(e)=\sum_{e \, \text{ends at} \, v}
f(e)
$$
at all vertices $v$ of $G$. Let $f(v)$ denote this quantity and let $\F(G)$ 
denote the {\em set of flows} of a digraph $G$.  
\end{definition}

If $\beta$ is a weight function on
the set of edges of $G$ and $f$ is a flow on $G$, then
\begin{itemize}
\item
the {\em weight} $\b(f)$ of $f$ is given by
$\b(f)=\prod_e \b(e)^{f(e)}$, where 
$\b(e)$ is the weight of the edge $e$. 
\item
The {\em multiplicity} at a vertex $v$ with outgoing edges $e_1,e_2,\dots$ is 
given by $\mult_v(f)=\binom{f(e_1)+f(e_2)+
\dots}{f(e_1),f(e_2),\dots}$, and the {\em multiplicity} of $f$ is given by
$\mult(f)=\prod_v \mult_v(f)$. 
\item
If $A$ is a subset of edges then we let $f(A)=\sum_{e\in A} f(e)$.
\end{itemize}

Let us summarize Foata-Zeilberger's theorem \cite[Theorem 1.1]{FZ} here. For 
the sake of completeness we include its proof in Appendix \ref{sec.FZ}.

\begin{theorem}
\lbl{thm.FZ}
If $(G,B)$ is a weighted digraph, then
\begin{eqnarray}
\lbl{eq.zeta1}
\zeta(G,B) 
&=& \frac{1}{\det(I-B)} \\
\lbl{eq.zeta2}
&=& \sum_{f \in \F(G)} \b(f)\, \mult(f) .
\end{eqnarray}
\end{theorem}

\begin{remark}
\lbl{rem.comb}
For $r=1$, the above Theorem states that
$$
\frac{1}{1-x}=\sum_{n=0}^\infty x^n
$$
where $x=b_{11}$. Thus, Theorem \ref{thm.FZ} is a version of the geometric
series summation.
\end{remark}

Another formula for the inverse of a determinant, the MacMahon Master Theorem, 
has been mentioned 
in the introduction. We will need its quantum version, proved in \cite{GLZ}.

In $r$-dimensional {\em quantum algebra} we have $r$ indeterminates 
$x_i$ ($1\leq i\leq r$),
satisfying the commutation relations $x_jx_i=qx_ix_j$ for all 
$1\leq i<j\leq r$. 
Further we are given a right-quantum matrix $A$. We assume
that the indeterminates of $A$ commute with the $x_i$'s. The following 
theorem has been proven
recently in \cite{GLZ}. 

\begin{theorem}
\lbl{thm.qMM}
Let $A$ be a right-quantum matrix of size $r$. For $1\leq i\leq r$ let 
$X_i=\sum_{j=1}^ra_{ij}x_j$,
and for any vector $(m_1,\dots, m_r)$ of non-negative integers let 
$G_A(m_1,\dots, m_r)$ be
the coefficient of $x_1^{m_1}x_2^{m_2}\dots x_r^{m_r}$ in 
$\prod_{i=1}^rX_i^{m_i}$. Then
$$
\sum_{m_1,\dots, m_r=0}^{\infty}G_A(m_1,\dots, m_r)=1/\Ferm(A).
$$
\end{theorem}

\section{The arc-graph of a knot projection}
\lbl{sec.arcgraph}

Given a knot projection $K$, we define the arc-graph $G_{\K}$ as follows:
\begin{itemize}
\item
The vertices of $G_{\K}$ are in 1-1 correspondence with the arcs of $\K$.
\item
The edges of $G_{\K}$ are in 1-1 correspondence with {\em transitions}
of $\K$, when we walk along $\K$ and we either go under a crossing 
(blue edges), 
or jump up at a crossing (red edges) according to Figure \ref{planarG2K}.
\end{itemize}

\begin{figure}[htpb]
$$ 
\psdraw{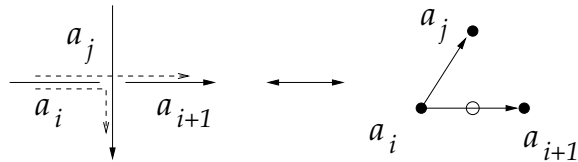}{3in}
$$
\caption{From a planar projection to the arc-graph. Transitions in 
the planar projection are indicated by dashed paths, and the corresponding 
edges in the arc-graph are blue (depicted with a small circle on them)
or red}\lbl{planarG2K}
\end{figure}

More formally, $G_{\K}$ is a weighted digraph defined as follows.

\begin{definition}
\lbl{def.arcgraph}
The {\em arc-graph} $G_{\K}$ has $r+1$ vertices $1,\dots, r+1$, $r+1$ {\em blue
directed edges} $(v,v+1)$ ($v$ taken modulo $r+1$) and $r+1$ {\em red 
directed edges} $(u,v)$, where at the crossing $u$ the arc that crosses over 
is labeled by $a_v$.

The vertices of $G_{\K}$ are equipped with a sign, where $\sign(v)$ is the
sign of the corresponding crossing $v$ of $\K$, and the edges of $G_{\K}$ are
equipped with a weight. The edge-weights are specified by matrix 
$W_{\K}= (\b_{ij})$
where 
$$
\b_{e}=
\begin{cases}
t^{-\sign(v)} & \text{if} \quad e=(v,v+1) \\
1-t^{-\sign(v)} & \text{if} \quad e=(v,u)
\end{cases}
$$
Here $t$ is a variable. Let $W_K$ denote the matrix obtained from $W_{\K}$
by deleting the last row and column.
Notice that $W_{\K}$ is formally stochastic
(i.e., the sum of the rows of $I-W_{\K}$ is zero), but $W_K$ is not.

Let $(G_K,W_K)$ denote the weighted digraph obtained by deleting the $r+1$ 
vertex
from $G_{\K}$, together with all edges to and from it. We let $V_K$ and $E_K$
denote the set of vertices and edges of $G_K$.
\end{definition}

It is clear from the definition that from every vertex of $G_{\K}$, the blue
outdegree is $1$, the red outdegree is $1$, and the blue indegree is $1$.
It is also clear that $G_{\K}$ has a Hamiltonian cycle that consists of all
the blue edges. We denote by $e_i^b$ ($e_i^r$) the blue (red) edge {\em leaving
vertex} $i$.

\begin{example}
\lbl{ex.1}
For the figure 8 knot we have:
$$
\psdraw{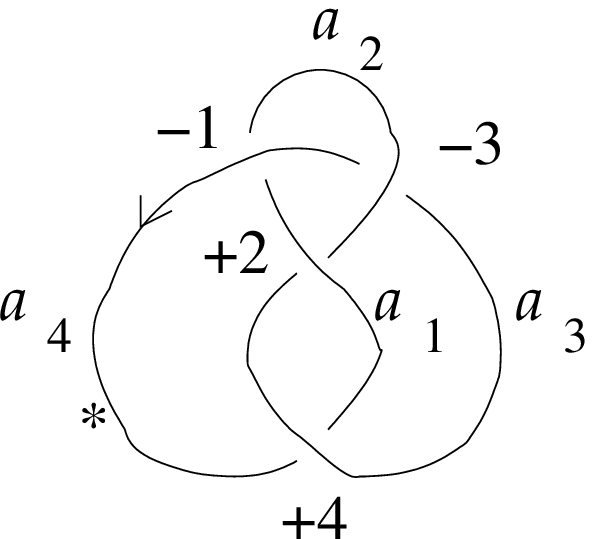}{1.5in}
$$
Its arc-graph $G_{\K}$ with the ordering and signs
of its vertices and  $G_K=G_{\K}-\{4\}$ are given by 
$$
\hspace{-1cm}
G_{\K}=
\psdraw{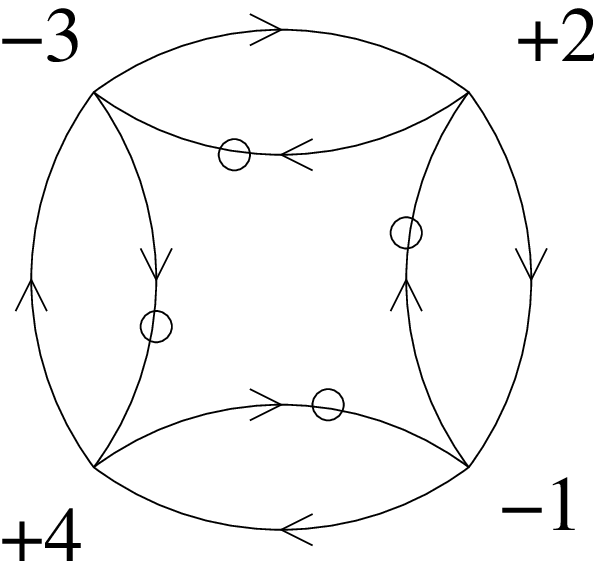}{1.3in} \qquad G_K=
\psdraw{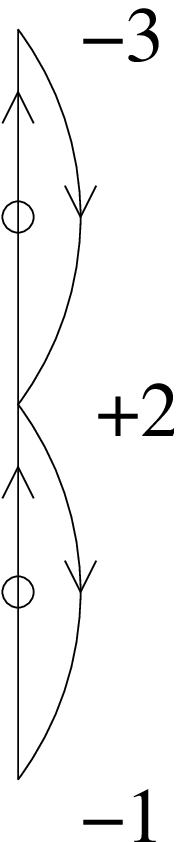}{0.3in}
$$  
where the blue edges are the ones with circles on them. Moreover, 
$$
W_{\K}=
\left[
\begin{matrix}
 0 & t & 0   & 1-t \\
 1-\bt & 0     & \bt & 0   \\
 0 & 1-t   & 0   & t \\
 \bt & 0   & 1-\bt   & 0   
\end{matrix}
\right] .
$$

\end{example}

\begin{definition}
\lbl{def.jj}
Let $\K$ be a knot projection. The {\em writhe} of $\K$, $\w(\K)$, is the 
sum of the signs of the
crossings of $\K$, and $\rot(K)$ is the {\em rotation number} of $K$, 
defined as follows: {\em smoothen} all crossings of $\K$, and consider the 
oriented circles that appear; one of them is special, marked by $\star$. 
The number of 
circles different from the special one whose orientation agrees with the 
special one, minus the number of circles whose orientation is opposite to the
special one is defined to be $\rot(K)$. 
We further let $\d(K,n)= 1/2 (n^2\w(\K)+n\rot(K))$, and $\d(K)= \d(K,1)$.
\end{definition}

We remark that we define $\rot(e)$ for each edge $e$ of $G_{\K}$ in 
Definition \ref{def.jjj}.

\section{The enhanced arc-graph and the Jones polynomial}
\lbl{sub.enhanced}

In order to express the Jones polynomial as a function of the arc graph, we 
need to enhance the arc-graph as follows.

\begin{definition}
\lbl{def.rotten}
{\rm (a)}
We introduce a {\em linear order} $<_v$ on the set of
edges of $G_K$ terminating at vertex $v$ as follows. Recall that
$v$ corresponds to an arc $a_v$ of $K$. 
If we travel on $a_v$ along the orientation of $K$, we 'see' one by one the 
arcs corresponding to starting vertices of red arcs 
entering $v$: this gives the linear order of red arcs entering $v$. 
Finally there is at most one blue edge entering vertex $v$, and we make it 
less than all the red edges entering $v$. 
\newline
{\rm (b)}
If $f$ is a flow on $G_K$, we define the {\em rotation} and {\em excess}
number of $f$ by:
\begin{equation}
\lbl{eq.rotexc}
\rot(f)=\sum_{e\in E_K} f(e)\rot(e),
\qquad
\exc(f)=\sum_{v \in V_K} \sign(v) f(e^b_v) \sum_{e<_v e^r_v}f(e),
\qquad \d(f)=\exc(f)-\rot(f),  
\end{equation}
where $V_K$ and $E_K$ are the set of vertices and edges of $G_K$, and $\rot(e)$
is defined in Definition \ref{def.jjj}.
\end{definition}

Let $\Sev(G)$ denote the set of all subgraphs $C$ of $G$ such that each 
component
of $C$ is a directed cycle. Note that $\Sev(G)$ may be identified with a 
finite subset of $\F(G)$ since the characteristic function of $C$ is a flow.

The next theorem, due to Lin-Wang, expresses the Jones polynomial of
a knot projection $\K$ in terms of the enhanced arc-graph of $K$. 
For the sake of completeness we include
its proof in Appendix \ref{sec.LW}.

\begin{theorem}
\lbl{thm.arcjones}\cite{LW} 
For every knot projection $\K$ we have:
$$
J(\K, t)= t^{\d(K)} \sum_{c \in \Sev(G_K)} t^{\d(c)} \b(c).
$$
\end{theorem}

We now give a similar formula for the {\em colored Jones function} $J_n$ 
of a knot. 
We will normalize the colored Jones function so that it is the constant
sequence $\{1\}$ for the unknot, and $J_n$ is the quantum group invariant
of knots that corresponds to the $(n+1)$-dimensional irreducible representation
of $\mathfrak{sl}_2$.

Recall the operation of {\em cabling} $\K^{(n)}$ the knot projection 
$\K$ $n$ times.
 Recall that $a_1,\dots, a_{r+1}$ are the arcs of $\K$.
Each $a_l$ is in the cabling replaced by $n^2$ arcs $a^l_{i,j}$, 
$i,j=1,\dots, n$,
with the agreement that the 'long arcs' obtained by cabling arc $a_k$ will be
$a^k_{1j}$, $j=1,\dots, n$ and the 'small arcs' obtained by cabling of crossing
$k$ will be denoted by $a^k_{ij}$ for $i=2,\dots, n$ and $j=1,\dots, n$.
Note that all crossings which replace the original crossing $k$ have the
same sign, equal to the sign of the crossing $k$.
(see figure before Lemma \ref{lem.evenn}). 

We further let $K^{(n)}$ denote the link obtained from $\K^{(n)}$ by deletion
of the $n$ special long arcs $a^{r+1}_{1j}$, $j=1,\dots, n$.

\begin{theorem}
\lbl{thm.arcjonesn}
For every knot $\K$ and every $n \in \BN$, we have
$$
J_n(\K, t)=t^{\d(K, n)} 
\sum_{c \in \Sev(G_{K^{(n)}})}
t^{\d(c)} \b(c) 
$$
where $G_{K^{(n)}}$ is the arc graph of $K^{(n)}$.
\end{theorem}

\begin{proof}
Let $V_n$ denote the $(n+1)$-dimensional irreducible representation of the
quantum group $U_q(\mathfrak{sl}_2)$, and let $v_n$ denote a highest weight 
vector of $V_n$. Then, there is an inclusion $V_n \longto \otimes^n V_1$ that 
maps $v_n$ to a nonzero multiple of $\otimes^n v_1$.  

The result follows since cabling $\K$ corresponds to tensor product of 
representations and since $\w(\K^{(n)})=n^2 \w(\K)$ and $\rot(K^{(n)})=
n\, \rot(K)$. 
\end{proof}

For an integer $m$, we denote by 
$$
(m)_q=\frac{q^m-1}{q-1}
$$ 
the {\em quantum integer} $m$. This defines the 
{\em quantum factorial} and the {\em quantum binomial coefficients} by
$$
(m)_q!=(1)_q (2)_q \dots (m)_q
\qquad
\binom{m}{n}_q=\frac{(m)_q!}{(n)_q!(m-n)_q!}
$$
for natural numbers $m,n$ with $n \leq m$. We also define
$$
\mult_q(f)=\prod_v\binom{f(v)}{f(e_v^b)}_{q^{\sign(v)}}.
$$

\begin{theorem}
\lbl{thm.main}
For every knot projection $\K$ we have:
$$
J_n(\K,t) = t^{\d(K, n)}\sum_{f\in \F(G_K)}\mult_{t}(f)
t^{\d(f)} \prod_{v \in V_K} t^{-\sign(v)nf(e^b_v)} 
\prod_{e \text{red};t(e)=v}\prod_{j=0}^{f(e)-1} (1-t^{-\sign(s(e))
(n- j- \sum_{e'<_v e}f(e))}). 
$$
\end{theorem}

\begin{remark}
\lbl{rem.jo1}
It simply follows that the
contribution of a flow $f$ to the sum in Theorem \ref{thm.main} is non-zero 
only if $f(v)\leq n$ for each vertex $v$. Thus, in the above sum, only 
finitely many terms contribute. As a result, when $n=1$,
Theorems \ref{thm.main} and \ref{thm.arcjones} coincide. 
\end{remark}

\section{Proof of Theorem \ref{thm.nonc}}
\lbl{sec.prpr}

Theorem \ref{thm.main} is used in this section to prove the main 
Theorem \ref{thm.nonc}.
In the rest of the paper we then prove Theorem \ref{thm.main} from 
Theorem \ref{thm.arcjonesn}.

\subsection{Row and column arcs order}
\lbl{sub.rco}

Recall that we fix a generic planar projection $\K$ of an oriented knot with 
$r+1$ crossings.
We order the arcs of $\K$ so that they appear
in increasing order as we walk in the direction of the knot,
and we denote by $K$ the {\em long knot} obtained by breaking the
arc $a_{r+1}$. We also order the crossings of $K$
so that arc $a_i$ ends at the $i$th crossing, for $i=1,\dots,r$.

\begin{definition}
\lbl{def.perm}
\begin{enumerate}
\item
We define two permutations $S,T$ on the set of the arcs of $K$ as follows.
For arc $a_i$ of $K$ let $T(i)= T(i,1),\dots, T(i,k_i)$ 
( $S(i)= S(i,1),\dots, S(i,k_i)$ respectively)
be the block of arcs of $K$ terminating at (starting from) $a_i$ and 
ordered along
the orientation of $a_i$. Let $T$ ($S$) be the permutation of the arcs of $K$
defined by $T= T(1,1),\dots, T(1,k_1),\dots, T(r,k_r)$ 
($S= S(r,k_r),\dots, S(r,1),\dots, S(1,1)$).
\item
We define permutation $R$ of the arcs of $K$ from $T$ as follows:
if $a_i$ appears in $S$ {\it before} the block $S(i)$ then 
replace $T(i,1),\dots, T(i,k_i)$ by $T(i,k_i),\dots, T(i,1)$.
\item
Similarly we define permutation $C$ of the arcs of $K$ from $S$ as follows:
if $a_i$ appears in $S$ {\it after} the block $S(i)$ then 
replace $S(i,k_i),\dots, S(i,1)$ by $S(i,1),\dots, S(i,k_i)$.
\end{enumerate}
\end{definition}

\begin{definition}
\lbl{def.b'}
We define matrix $B'_K= (\gamma_{ij})$ to be obtained from $B_K$ by 
taking the rows in the $R$ order
and the columns in the $C$ order.
\end{definition}

We consider the commutation relations between the variables appearing in 
$B'_K$ as 
in the Definition \ref{def.com}. In particular, $B'_K$ is right-quantum.


\subsection{Flows on $G_K$ and monomials of $G_{B'_K}(m_1,\dots m_r)$}
\lbl{sub.mon}

We interpret each entry $\gamma_{ij}$ with no $z$ indeterminate as arc 
$(ij)$ of the arc-graph $G_K$.
Then each monomial in $G_{B'_K}(m_1,\dots, m_r)$ corresponds to a flow on 
$G_K$
with $indeg(i)= outdeg(i) = m_i, i=1,\dots, r$. If $f$ is such a flow, we 
denote by $G(f)$
the sum of all monomials of $\sum G_{B'_K}(m_1,\dots, m_r)$ corresponding to 
$f$. Summarizing we can write
 
\begin{observation}
\lbl{obs.fl}
$$
\sum_{m_1,\dots, m_r=0}^{\infty} G_{B'_K}(m_1,\dots, m_r)= \sum_f G(f).
$$
\end{observation}

We denote by $C(f)$ the canonical monomial of a product 
(in arbitrary order) of the entries of $B'_K$
corresponding to the edges of $G_K$, where the entry corresponding to 
each edge $e$
appears $f(e)$ times.

\begin{observation}
\lbl{obs.hr}
Let $C$ be a summand of $\prod_{i=1}^r(\sum_{j}\gamma_{ij}x_j)^{m_i}$, which 
contains $m_i$ indeterminates $x_i, i=1,\dots, r$ and contributes to $G(f)$.
For $1\leq v\leq r$ and $1\leq j\leq f(e^r_v)$ let $c(C,v,j)$ be the 
number of $\gamma_{e^b_v}$'s 
which need to be commuted through the $j$-th occurrence of $\gamma_{e^r_v}$ 
in order to get
$C(f)$ from $C$.
Then
$$
C= x_1^{m_1}x_2^{m_2}\dots x_r^{m_r}q^{-\rot(f)}C(f)q^{\exc(f)}
\prod_{v=1}^r\prod_{j=1}^{f(e^r_v)}
q^{\sign(v)c(C,v,j)}.
$$
\end{observation}
\begin{proof}
Let $X_{ij}= \gamma_{ij}x_j$. Hence $C$ is a summand of the coefficient of 
$x_1^{m_1}x_2^{m_2}\dots x_r^{m_r}$ in $\prod_{i=1}^r(\sum_{j}X_{ij})^{m_i}$. 
For each $j$ fixed the $\gamma_{ij}$'s appear ordered in $C$.
In each canonical monomial, the $\gamma_{ij}$'s appear ordered by the 
second coordinate, and then by the 
first coordinate. Hence, in order to get a canonical monomial times 
$x_1^{m_1}x_2^{m_2}\dots x_r^{m_r}$ from a summand of 
$\prod_{i=1}^r(\sum_{j}X_{ij})^{m_i}$,
we only need to commute $X_{ij}$'s so that they are ordered by the 
second coordinate.
This means: if $a_i$ appears in $S$ {\it before} ({\it after} respectively)
the block $S(i)$ 
then $a_i$ appears in $C$ {\it before} ({\it after} respectively )the 
block $S(i)$ 
BUT $a_{i-1}$ appears in $R$ {\it after} ({\it before} respectively) 
the block $T(i)$. Hence 
we need to commute 
\begin{enumerate}
\item[1.]
Each $X_{i-1,i}$ through each $X_{j-1,j}, j\in S(i)$ and
each $X_{k-1,i}$ through each  $X_{j-1,j}, j,k\in S(i), R(k)>R(j)$.
The commutation in $B'_K$ is such that we acquire each time $q^{\sign(j-1)}$.
Hence we acquire in total $q^{\exc(f)}$ since we recall that
$\exc(f)=\sum_{v} \sign(v) f(e^b_v) \sum_{e< e^r_v}f(e)$.
\item[2.]
Each $X_{k-1,i}$ through each  $X_{k-1,k}, k\in S(i)$. The commutation in 
$B'_K$ is such that 
we acquire  in total
$q^{\sign(k-1)(c(C,k-1,1)+\dots c(C,k-1,f(e^r_{k-1}))}$. 
\item[3.]
 The commutation in $B'_K$ is such that if $i<i', j<j'$ and $X_{i,j'},
X_{i',j}$ do not
appear in one of the previous two cases then $X_{i,j'},X_{i',j}= X_{i',j},
X_{i,j'}$.
\end{enumerate}
This finishes the proof.
\end{proof}

\begin{corollary}
\lbl{cor.h}
$$
G(f)= C(f)q^{\d(f)}\prod_{v=1}^r
\sum_{f(e^b_v)\geq c_1\geq\dots \geq c_{f(e^r_v)}\geq 0}
q^{\sign(v)(c_1+\dots c_{f(e^r_v)})}.
$$
\end{corollary}

Since
$$
\sum_{f(e^b_v)\geq c_1\geq\dots \geq c_{f(e^r_v)}\geq 0}
q^{\sign(v)(c_1+\dots c_{f(e^r_v)})}=
\binom{f(v)}{f(e^b_v)}_{q^{\sign(v)}},
$$
we have 

\begin{corollary}
\lbl{cor.skoro}
$$
G(f)= q^{\d(f)}\mult_q(f)C(f).
$$
\end{corollary}

\begin{proof}(of Theorem \ref{thm.nonc})
Theorem \ref{thm.main} tells us that 
$$
J_n(K,q) = q^{\d(K,n)}\sum_{f\in \F(G_K)}\mult_{t}(f)
t^{\d(f)} \prod_{v \in V_K} t^{-\sign(v)nf(e^b_v)}
\prod_{e\text{red};t(e)=v}\prod_{j=0}^{f(e)-1} 
(1-t^{-\sign(s(e))(n- (\sum_{e'<_v e}f(e))-j)}).  
$$
Comparing this with the definition \ref{def.eval} of the n-evaluation and 
using 
Theorem \ref{thm.qMM}, Observation \ref{obs.hr} and Corollary \ref{cor.skoro},
we can see that Theorem \ref{thm.nonc} follows.
\end{proof}

\section{Cabling of the arc graph}
\lbl{sec.stateCJ}

Recall the operation of {\em cabling} $\K^{(n)}$ the knot projection $\K$ $n$ 
times.
 Recall that $a_1,\dots, a_{r+1}$ are the arcs of $\K$.
Each $a_l$ is in the cabling replaced by $n^2$ arcs $a^l_{i,j}$, 
$i,j=1,\dots, n$,
with the agreement that the 'long arcs' obtained by cabling arc $a_k$ will be
$a^k_{1j}$, $j=1,\dots, n$ and the 'small arcs' obtained by cabling of crossing
$k$ will be denoted by $a^k_{ij}$ for $i=2,\dots, n$ and $j=1,\dots, n$.
Note that all crossings which replace the original crossing $k$ have the
same sign, equal to the sign of the crossing $k$.
We make the following agreement: assume the parallel arcs 
$a^l_{i,1},\dots, a^l_{i,n}$
go horizontally from left to right. Then $a^l_{i,1}$ is the upmost one.
See figure below for part of Example \ref{ex.1} and of its
3-cabled version $\K^{(3)}$:
$$
\psdraw{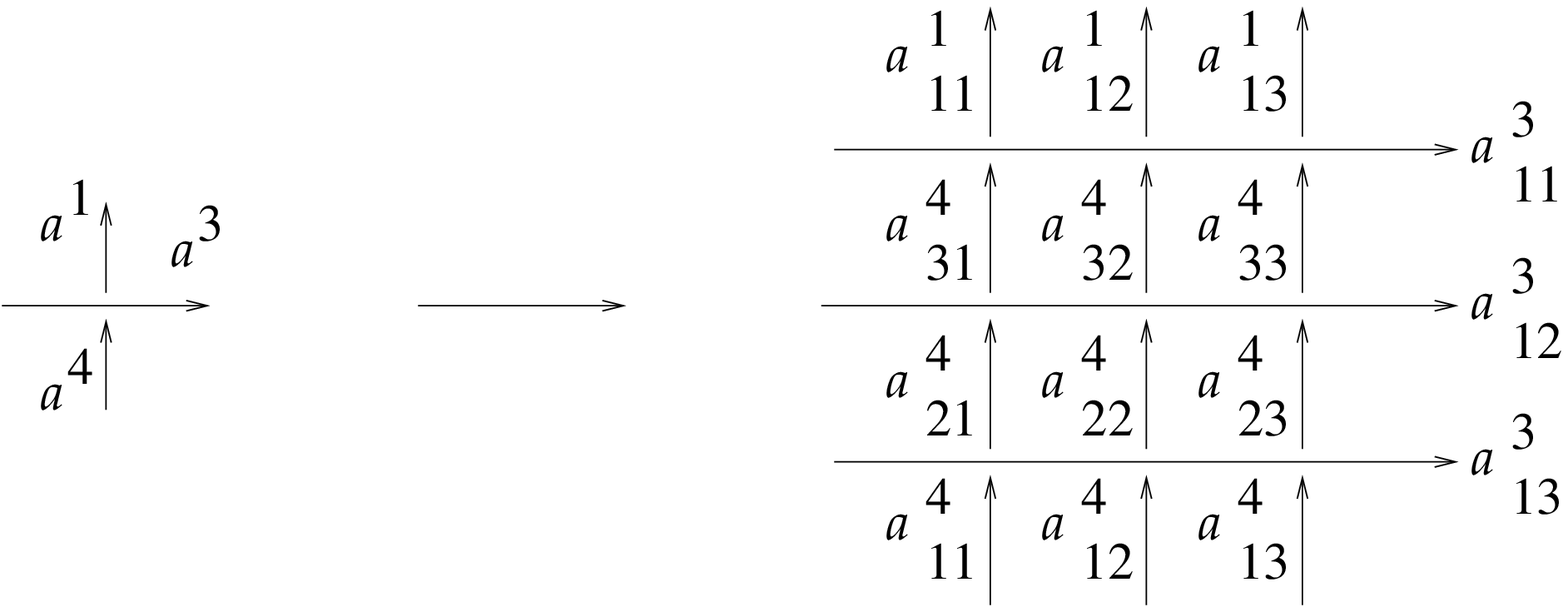}{4in}
$$

We further let $K^{(n)}$ denote the link obtained from $\K^{(n)}$ by deletion
of the $n$ special long arcs $a^{r+1}_{1j}$, $j=1,\dots, n$.

Next we consider the arc-graph of the cabling of $\K$.
For example, part of the red-blue digraph $G_{\K}$ of Example \ref{ex.1} and 
of its 3-cabled version is depicted as follows:
$$
\psdraw{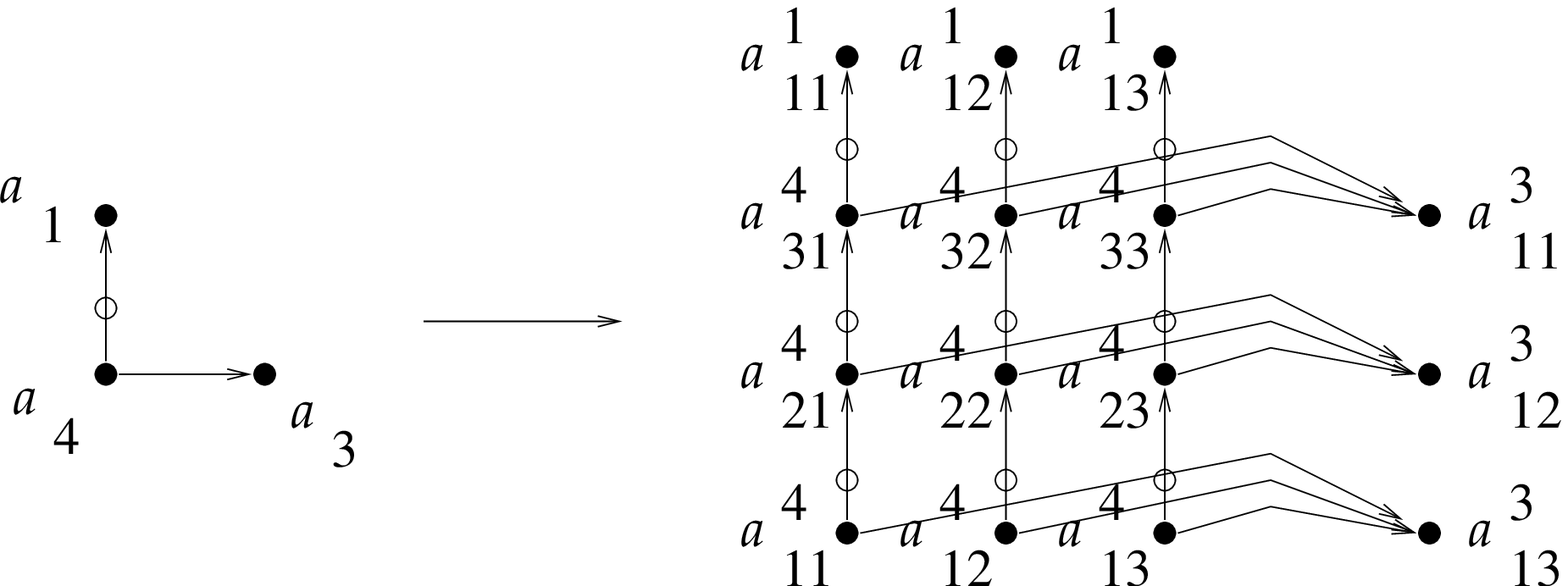}{3.5in}
$$
where vertical edges are blue and horizontal edges are red.

We now define an n-cabling $G_K^{(n)}$ of the arc-graph $G_K$. 
Cabling of a planar projection is a local operation, and
so is cabling of a digraph. In the language of combinatorics, we blow up
the vertices of $G$ using a suitable {\em gadget}. For a similar discussion,
see also \cite[Section 4]{GL}.

\begin{definition}
\lbl{def.cabled2}
Fix a red-blue arc-graph $G_K$. Let $G_K^{(n)}$
denote the digraph with vertices $a^v_j$ for $v$ a vertex of $G_K$ 
and $j=1,\dots,n$.
$G_K^{(n)}$ contains blue directed edges $(a^l_j,a^{l+1}_j)$ with weight
$t^{-\e n}$ (where $\e \in \{-1,+1\}$ is the sign of the crossing $l$)
for each $l=1,\dots,r-1$ and $j=1,\dots,n$. Moreover,
if $(a_k,a_l)$ is a red directed edge of $G_K$, then $G_K^{(n)}$ contains
red edges $(a^k_i,a^l_j)$ for all $i,j=1,\dots n$ with weight
$t^{(j-1)}(1-t)$ resp. $t^{-(n-j)}(1-t^{-1})$, if the sign of the $i$ 
crossing is $-1$ resp. $+1$. Notice that the weights of the red edges are 
independent of the index $i$. 
\end{definition}

\begin{lemma}
\lbl{lem.evenn}
There is a 1-1 correspondence
$$
\{\text{\rm{admissible subgraphs of}} \,\,
G_{K^{(n)}}\} \longleftrightarrow 
\{\text{\rm{admissible subgraphs of}} \,\, G_K^{(n)} \}.
$$
\end{lemma}

We will denote the set of admissible even subgraphs of $G_K^{(n)}$ by 
$\Sev_n(G_K)$.

\begin{proof}
Denote by $p^k_j$ path $(a^k_{1j},a^k_{2j}, \dots, a^k_{nj})$ of $n-1$
blue edges in $G_{K^{(n)}}$, $k=1,\dots, r$ and $j=1,\dots,n$. 
There is a natural map $G_{K^{(n)}}\longto G_K^{(n)}$ which contracts each 
directed path $p^k_j$ into its initial vertex, and deletes all vertices 
$a^{r+1}_{ij}$.
Forgetting the weights, it is clear that the result of the contraction
coincides with $G_K^{(n)}$.

$G_{K^{(n)}}$ has two types of vertices: $a^k_{ij}$ for
$k=1,\dots, r+1$ and $i,j=1,\dots,n$ and $i \neq 1$ (call these white)
and $a^k_{1j}$, $k< r+1$ (call these black). The indegree of a white vertex 
is $1$, but the indegree of a black vertex may be higher. The black vertices
are the initial vertices of the paths $p^k_j$, hence the vertices
of $G_K^{(n)}$. Let $E^{\{n\}}$ and $E^{(n)}$ denote the set of edges of 
$G_{K^{(n)}}$ 
and $G_K^{(n)}$ respectively. Then each edge $e$ of $E^{(n)}$ replaces the 
unique 
directed path $P_e$ of $G_{K^{(n)}}$ between the corresponding black 
end-vertices of $e$, 
which contains no other black vertices.
If $E \subset E^{\{n\}}$ is an admissible even subgraph of $G_{K^{(n)}}$ 
then $E$ is a vertex-disjoint union of directed cycles of  $G_{K^{(n)}}$
and each directed cycle may be decomposed into directed paths between the 
black vertices.
If each such directed path is replaced by a directed edge, we get an admissible
even subgraph $E'$ of $G_K^{(n)}$. This gives the 1-1 correspondence between 
the admissible even subgraphs without the weights. To realize that the 
weights are correct as well, we only need to compare the product of the 
weights in $G_{K^{(n)}}$ along $P_e$ with the weights of $e$ in $G_K^{(n)}$.
\end{proof}

Theorem \ref{thm.arcjonesn} and Lemmas \ref{lem.arc}, \ref{lem.evenn} 
imply that:

\begin{lemma}
\lbl{lem.arcjonesnn}
For every knot $\K$ and every $n \in \BN$, we have
$$
J_n(\K, t)=
t^{n \d(K)} \sum_{c \in \Sev_n(G_K)} t^{\d(c)} \b(c). 
$$
\end{lemma}

Our next task is to figure out $\d(c)=\exc(c)-\rot(c)$ for $c\in \Sev_n(G_K)$. 

The following lemma is clear from the Definition \ref{def.rotten}: 

\begin{lemma}
\lbl{lem.rot}
If $f$ is a flow on $G_K$ and $\ti f$ is a lift
of $f$ to flow on $G_K^{(n)}$,
for some $n$, then $\rot(f)=\rot(\ti f)$.
\end{lemma}

\subsection{Comparison of excess numbers}
\lbl{sub.coexcess}

Given an admissible subgraph $c$ in $G_K^{(n)}$, let $f$ be the 
corresponding flow in $G_K$, to which $c$ projects, under the projection
$$
\pi: G_K^{(n)} \to G_K.
$$ 
In this section, we 
compare $\exc(f)$ (in Definition \ref{def.rotten}) with
$\exc(c)$.

As we will see, the two excess numbers do not agree. In this section
we will determine their difference.

We begin by introducing a partial ordering $\prec$ on the set of edges
of $G_K^{(n)}$. We warn the reader that this ordering is different 
from the ordering $<_v$ of the edges of $G_K$ entering vertex $v$,
introduced in Definition \ref{def.rotten}.

\begin{definition}
\lbl{def.orderf}
Consider two edges $e$ and $e'$ of $G_K^{(n)}$ which start at the vertices
$a^i_j$ and $a^{i'}_{j'}$  of $G_K^{(n)}$. We say that
$e \prec e'$ if
\begin{itemize}
\item
$e,e'$ end at the same vertex $v$ and $\pi(e) <_v \pi(e')$ in $G_K$, or
\item
$i=i'$ and $\sign(i)=+$ and $j < j'$, or
\item
$i=i'$ and $\sign(i)=-$ and $j' < j$
\end{itemize}
\end{definition}

Recall that $c\in \Sev_n(G_K)$ ($c$ admissible) if and only if
 $c$ is a collection of vertex disjoint directed cycles of $G_K^{(n)}$.
Hence the ordering on the edges of $c$ defined in Definition \ref{def.orderf}
induces a total ordering on each $\pi^{-1}(e)$, $e$  edge of $G_K$.

This total ordering may be seen from the cabling of the knot in the same way 
as the ordering $<_v$ of Definition \ref{def.rotten} may be seen from the knot:
if we travel along an arc $a^i_j$, we see one by one the arcs corresponding
to the starting vertices of edges of $\pi^{-1}(e)$, where $e$ is an edge of 
$G_K$.
This agrees with the total ordering on $\pi^{-1}(e)$ induced by $\prec$;
see Figure before \ref{lem.evenn}.

\begin{definition}
\lbl{def.XY}
Consider two edges $e$ and $e'$ of $G_K^{(n)}$ which end at the vertices
$a^i_j$ and $a^{i'}_{j'}$  of $G_K^{(n)}$. 
$$
X(e,e')=
\begin{cases}
1 & \text{if} \quad e, e'= \text{ red}, \quad 
e' \prec e, \quad \sign(s(e))= +, \quad j < j' \\
1 & \text{if} \quad e=\text{red}, \,\, e'=\text{blue},\,\, \\
& \pi(e), \pi(e') \,\, \text{do not start at the same vertex}, \\ 
& e' \prec e, \quad \sign(s(e))= +, \quad j < j' \\
1 & \text{if} \quad e,e'= \text{ red}, \quad
 e' \prec e, \quad \sign(s(e))= -, \quad 
j' <  j \\
1 & \text{if} \quad e=\text{red}, \,\, e'=\text{blue},\,\, \\
& \pi(e), \pi(e') \,\, \text{do not start at the same vertex}, \\ 
&  e' \prec e, \quad \sign(s(e))= -, \quad 
j' <  j \\
0 & \text{otherwise}
\end{cases}
$$
$$
Y(e,e')=
\begin{cases}
1 & \text{if} \quad e=\text{red}, \,\, e'=\text{blue},\,\, \\
& \pi(e), \pi(e') \,\, \text{start at the same vertex}, \quad e \prec e' \\
0 & \text{otherwise}
\end{cases}
$$
\end{definition}

\begin{lemma}
\lbl{lem.exc2}
Let $c$ be an admissible subgraph of $G_K^{(n)}$. 
Denote by $f$ the flow on $G_K$
which is the projection of $c$ to $G_K$. Then 
$$
\exc(c)=\exc(f) + \sum_{e}\sign(s(e)) \left(\sum_{e'}
X(e,e')+Y(e,e') \right) 
$$
where the summations of $e$ and $e'$ are over the set of edges of $c$ and
$s(e)$ denotes the starting vertex of $e$. 
\end{lemma}

\begin{proof}
Consider a crossing $v$ of $K$, and the corresponding $n^2$ crossings
of $K^n$. 
We count the contribution to  
$\exc(c)$ of pairs $(e,e')$ of edges of $c$ such that
\begin{itemize}
\item
$e$ projects to $e^b_v$ (the blue edge that starts at $v$),
and $e'$ does not project 
 to $e^r_v$ (the red edge that starts at $v$). This
gives $\exc(f)$.
\item
$e$ projects to $e^r_v$, and $e'$ does not project to $e^b_v$. This
gives the $X$-term in the formula. 
\item
$e$ projects to $e^r_v$, and $e'$  projects to $e^b_v$. 
This gives the $Y$-term in the formula. 
\end{itemize}
\end{proof}

\section{Sortings and multiplicities of flows}
\lbl{sec.categorifying}

\subsection{Sortings}
\lbl{sub.CAlex}

In this section we introduce one of our key tools, which is a categorification
of multiplicities of the flows on  $G_K$. Let $f$ be a flow on $G_K$. If $e$ 
is an edge of $G_K$
then we let $F(e)\subset F$ be the set of $f(e)$ copies of $e$; we choose
an {\em arbitrary} total order on each $F(e)$.

Let $F= \cup_{e\in E_K}F(e)$ and let $F_r\subset F$ consists of the union of 
$F(e)$, $e$ red.
Further let $F^+_r$ denote the subset of $F_r$ consisting of the red edges 
which leave a vertex with 
$+$ sign, and we let $f^+_r=|F^+_r|$. Analogously we define  $F^-_r, \dots $.

\begin{definition}
\lbl{def.conf}
Fix a flow $f$ on $G_K$. A {\em sorting} $C$ of $f$
is a function
$$
C : \text{Vertices}(G_K) \to 2^{F_r}
$$
such that 
\begin{itemize}
\item
$C_1$ is a collection of red edges that terminate in vertex $1$, of
cardinality $f(e^b_1)$.
\item
For each $2\leq i \leq r$,
$C_{i} \sub C_{i-1}\cup \{e\in F_r; e$ terminates in vertex $i\}$ 
of $f(e^b_i)$ elements.
\end{itemize}
Let $\C(f)$ denote the set of all sortings of $f$. 
\end{definition}

\begin{lemma}
\lbl{lem.conf1}
Every flow $f$ has $\mult(f)$ sortings.
\end{lemma}

\begin{proof}
Use that $\mult(f,r)=1$, $\{e\in F_r; e$ terminates in vertex $1\} 
=\{e\in F; e$ 
terminates in vertex $1\}$ and for each $2\leq i <r$, $\sum f(e): e$ terminates
at vertex $i$ equals $|C_{i-1}\cup \{e\in F_r; e$ terminates in vertex $i\}|$.
\end{proof}

\begin{definition}
\lbl{def.Ifn}
We define  $\I(f,n)=\{0,\dots,n-1\}^{F_r}$. 
If $v\in \I(f,n)$ then we define
$f^-_r(v)=\sum_{e\in F^-_r}v_e$ and we define $f^+_r(v)$ analogously. 
\end{definition}

\begin{definition}
\lbl{def.filled}
{\rm (a)}
Fix a flow $f$ on $G_K$ and a natural number $n$.
An {\em n-sorting of $f$ } is a pair
$P=(C,v)$ where $C\in \C(f)$ and $v\in \I(f,n)$. 
\newline
{\rm (b)} If $P=(C,v)$ is an n-sorting of $f$ 
then we define its weight $b(P)$ to be
$$
b(P)= t^{n(f^-_w-f^+_w)}
 (1-t)^{f^-_r}t^{f^-_r(v)}(1-t^{-1})^{f^+_r}t^{-(f^+_r)(n-1)+ f^+_r(v)}.     
$$
{\rm (c)}
Let $\C_n(f)$ denote the {\em set of all n-sortings of $f$}.
\end{definition}

The following lemma states that 
the n-sortings of $f$ categorify 
multiplicities and weights of flows.

\begin{lemma}
\lbl{lem.geomseries}
 For every flow $f$ on $G_K$ and every $n$
we have
$$
 \b(f)|_{t\to t^n} \mult(f)=
\sum_{P\in \C_n(f)} b(P).
$$
\end{lemma}

\begin{proof}
It follows by Lemma \ref{lem.conf1} that
$$
\sum_{P\in \CI(f,n)} b(P)= \mult(f) t^{n(f^-_b-f^+_b)}(1-t)^{f^-_r}
(1-t^{-1})^{f^+_r}
 \left(\sum_{v\in \I(f,n)}t^{f^-_r(v)}t^{-(f^+_r)(n-1)+ f^+_r(v)}\right)=
$$
by a simple rearrangement
$$
\mult(f) t^{n(f^-_b-f^+_b)}
 \left((1-t)^{f^-_r}\sum_{v\in \I(f,n)}t^{f^-_r(v)}\right)
\left((1-t^{-1})^{f^+_r}\sum_{v\in \I(f,n)}t^{-(f^+_r)(n-1)+ f^+_r(v)}\right)=
$$
by a {\em geometric series summation}
$$
 \b(f)|_{t\to t^n} \mult(f).
$$
\end{proof}

\subsection{Sortings and the Jones polynomial}
\lbl{sub.CJones}

Here we define admissible sortings and give a formula for the colored
Jones function in terms of them.

\begin{definition}
\lbl{def.AConf}
Fix a flow $f$ of $G_K$ and a natural number $n$. Let 
$P=(C,v)$ be an n-sorting of $f$. We say that $P$ is {\em admissible} if
\begin{itemize}
\item
For every two edges $e,e' \in F_r$ such that $v_e=v_{e'}$ and 
$e$ ends in vertex $i$ and $e'$ ends in vertex $j$ and $j \geq i$,
there exists an $l$, $i \leq l < j$ such that $e \not\in C_l$.
\end{itemize}
We denote by $\AC_n(f)$ the set of all admissible n-sortings of $f$,
and by $\Sev_n(G_K, f)$ the set of all lifts of $f$ to admissible subgraphs 
of $G_K^{(n)}$.
\end{definition}

The next lemma explains the notion of admissible sortings.

\begin{lemma}
\lbl{lem.ACconf}
There is a bijection $\Phi$ from $\AC_n(f)$ to $\Sev_n(G_K, f)$
such that $\b(\Phi(P))= b(P)$.
Moreover, if $P\in \AC_n(f)$  and $\pi$ is the projection to  $G_K$ 
then for each $e\in E_K$, the fixed total order on $F(e)$ agrees 
with the total order $(\pi^{-1}(e), \prec)$ introduced in Definition 
\ref{def.orderf}.
\end{lemma}

\begin{proof}
Let $P=(C,v)$, $C=(C_1,\dots, C_{r})$, and $P\in \AC_n(f)$.
In order to define
$\Phi(P)$ we will define the image $\Phi(P,e)$ for each $e\in F$.

First we  determine the ends of  the lifts of the red edges
as follows: if $e\in F_r$ then we let $t(\Phi(P,e))= a^{t(e)}_{v(e)+1}$.

Next we determine the blue edges of $\Phi(P)$ as follows:
if $1\leq i \leq r$ then we let 
$$
\{s(\Phi(P,e))\, |\,\, e\in F(e^b_i)\}= \{a^i_{v_e+1}\, | \,\, e\in C_i\}.
$$
This determines the beginnings of the  blue edges,
and hence also the ends of the blue edges.

 It remains to specify the beginnings of the lifts of the red edges. 
Since $P$ is admissible, observe that for each $1\leq i\leq r$, 
there are exactly $f(e^r_i)$ vertices $a^i_j$ of indegree $1$ in current
$\Phi(P)$. Hence it remains to make each of them starting vertex of exactly
one edge $\Phi(P,e)$, $e\in F(e^i_r)$. This
is uniquely determined by the 'moreover' part of the Lemma. 
This finishes the definition of $\Phi$. The equality for the weights 
follows easily, and the moreover part of the Lemma directly from the 
definition of $\Phi$. To finish the proof we find the inverse to $\Phi$.

Let $c\in \Sev_n(G_K, f)$. We construct $\Phi^{-1}(c)=(C(c), v(c))$ as follows:
 Let $e$ be an edge of $G_K$. There is an order preserving bijection between 
the 
fixed total ordering $(F(e),<)$ and  $(\pi^{-1}(e), \prec)$. If $e'$ is an 
edge 
of $c$ then we let $e'_F$ be the corresponding edge of $F$.

First let $e'$ be a red edge of $c$. We let $v(c)_{e'_F}= j$ where 
$t(e')=a^i_{j+1}$.
Hence $v(c)$ encodes the ends of the red edges of $c$.

Next we  define a predecessor $p(e)$ for each edge $e$ of $c$. 
If $e$ red then $p(e)=e$. If $e$ blue then $p(e)$ is the red edge of $c$ 
which terminates in the starting vertex of the longest blue path of $c$ 
whose last edge is $e$. Note that $p(e)$ always exists and is unique
since $c$ is admissible.

Finally for $1\leq i \leq r$ let $C_i=\{p(e)_F; e$ edge of $c$ that 
terminates in 
some $a^i_j$ that is a starting vertex of a blue edge of $c \}$.
This finishes the construction of $\Phi^{-1}$. 
\end{proof}

\begin{theorem}
\lbl{thm.sic5}
We have:
$$
J_n(K)(t) = t^{\d(K,n)}
\sum_{f\in \F(G_K)}
t^{\d(f)}\sum_{P\in \AC_n(f)}t^{\exc(P)}b(P) ,
$$
where $\exc(P)=\exc(\Phi(P))-\exc(f)$. 
\end{theorem}

\begin{proof}
We have:
\begin{eqnarray*}
J_n(K)(t) &=& 
t^{\d(K,n)} \sum_{c \in \Sev_n(G_K)} t^{\d(c)} \b(c) \\
&=&
t^{\d(K,n)} \sum_{f\in \F(G_K)}t^{\d(f)}
\sum_{c \in \Sev_n(G_K,f)}   t^{\exc(c)-\exc(f)}\b(c) \\
&=&
t^{\d(K,n)}
\sum_{f\in \F(G_K)}t^{\d(f)} \sum_{P\in \AC_n(f)}t^{\exc(P)}b(P)
\end{eqnarray*}
\end{proof}

\subsection{Proof of Theorem \ref{thm.main}}
\lbl{sub.proof}

\begin{definition}
\lbl{def.prd}
Let $e\in F_r$. We define set $P(f,e)$ as follows:
if $e\in F(e'), e_1\in F(e'_1)$ then $e_1\in P(f,e)$ if $t(e')= t(e'_1)= v$
and $e'_1<_v e'$, or $e'= e'_1$ and $e_1< e$ in our fixed total order of 
$F(e)$.
\end{definition}

\begin{definition}
\lbl{def.okk}
Let $e \in F_r$. We define
\begin{itemize}
\item
$\defect_1(C,v,e)= |\{e'\in P(f,e): v_{e'}< v_{e} \}|$,
\item
$\defect_2(C,v,e)= |\{e'\in C_{d(e)}: v_{e'}< v_{e} \}|$,
where $d(e)$ is the biggest index such that $d(e)\geq t(e)$
and $e\not\in C_{d(e)}$.
\end{itemize}
Recall that $t(e)$ denotes the terminal vertex of $e$.
\end{definition}

\begin{proposition}
\lbl{prop.gk2}
Let $P= (C,v)$ be an n-sorting of $f$. Then
$$
\exc(P)=\sum_{e\in F_r}\d_1(e)+\d_2(e),
$$
where 
$$
\d_1(e)=\begin{cases}
|P(f,e)|-\defect_1(C,v,e) & \text{ if $e$ starts in a $+$ vertex}
\\
-\defect_1(C,v,e) & \text{ if $e$ starts in a $-$ vertex},
\end{cases}
$$
$$
\d_2(e)= 
\begin{cases}
|C_{d(e)}|- \defect_2(C,v,e) &\text{ if } \quad \sign(d(e))=+ \\
- \defect_2(C,v,e) & \text{ if} \quad \sign(d(e))=-.
\end{cases}
$$
\end{proposition}
\begin{proof}
Let $e\in F_r$ and first assume $\sign(s(e))= +$. Then 
$$
\d_1(e)= |P(f,e)|- \defect_1(C,v,e)= 
|\{e': e'\in P(f,e)\cap F_r, v_{e}< v_{e'} \}|+ 
|\{e' \in C_{t(e)-1}: v_{e}< v_{e'} \}|.
$$
This equals, by definition \ref{def.XY} of function $X$ and by definition of 
bijection $\Phi$ in 
lemma \ref{lem.ACconf},  
$$
\sign(s(\Phi(e)))\sum_{e'\in F}X(\Phi(e), \Phi(e')).
$$
We proceed analogously if $\sign(s(e))= -$. Hence
$$
\sum_{e\in F_r}\d_1(e)= \sum_{e\in F_r}\sign(s(\Phi(e)))
\sum_{e'\in F}X(\Phi(e), \Phi(e')).
$$
Next we denote, for $e\in F_r$, by $D(e))$ the red edge of $\Phi(P)$ that 
starts at vertex
$a^{d(e)}_{v(e)+1}$.
Note that $D$ is a bijection between $F_r$ and the set of the red edges of 
$\Phi(P)$. 

Now let $\sign(d(e))=+$. Then $\d_2(e)$ equals the number of $e'\in C_{d(e)}$ 
such that
$v_e< v_{e'}$. This equals $\sign(d(e))\sum_{e'\in F}Y(D(e), \Phi(e'))$.
Again the case $\sign(d(e))= -$ is analogous.

Hence we get 
$$
\sum_{e\in F_r}\d_2(e)= \sum_{e\in F_r}\sign(d(e))\sum_{e'\in F}Y(D(e), 
\Phi(e')).
$$
This finishes the proof by lemma \ref{lem.exc2}.

\end{proof}

\begin{proof}(of Theorem \ref{thm.main})

We use Theorem \ref{thm.sic5}, Lemma \ref{lem.ACconf} and 
Proposition \ref{prop.gk2}:

$$
J_n(K)(t) = t^{\d(K,n)}
\sum_{f\in \F(G_K)}
t^{\d(f)}\sum_{P\in \AC(f,n)}t^{\exc(P)}b(P)=
$$
$$
t^{\d(K,n)}\sum_{f\in \F(G_K)}t^{\d(f)}
t^{n(f^-_b-f^+_b)}(1-t)^{f^-_r}(1-t^{-1})^{f^+_r}
$$
$$
 \prod_{e\in F^+_r} t^{-(n-1-|P(f,e)|)}
\prod_{e\in F_r, \sign(d(e))=+}t^{|C_{d(e)}|}
\sum_{(C,v) \in \AC(f,n)}\prod_{e\in F_r} t^{v_e-\defect_1(C,v,e)-
\defect_2(C,v,e)}.
$$



Let us recall that
$$
\mult(f)_{q}= \prod_{v=1}^{r-2}\binom{f(v)}{f(e^b_v)}_{q^{\sign(v)}}.
$$

At this point, we will use Appendix \ref{sec.count}. By Theorem 
\ref{thm.sorstr} and 
Theorem \ref{thm.str10} we get 
$$
J_n(K)(t) = t^{\d(K,n)}
\sum_{f\in \F(G_K)}t^{\d(f)}
t^{n(f^-_b-f^+_b)}(1-t)^{f^-_r}(1-t^{-1})^{f^+_r}  t^{-(n-1)f^+_r} 
$$
$$
\prod_{e\in F^+_r} t^{|P(f,e)|}\prod_{e\in F_r, \sign(d(e))=+}t^{|C_{d(e)}|}
 \prod_{v=1}^{r}\binom{f(v)}{f(e^b_v)}_{t^{-1}} 
\prod_{e\in F_r}(n-|P(f,e)|)_t =
$$
$$
t^{\d(K,n)}\sum_{f\in \F(G_K)}t^{\d(f)}\prod_{v=1}^{r}
\binom{f(v)}{f(e^b_v)}_{t^{-1}}
\prod_{v:\sign(v)=+}t^{f(e^r_v)f(e^b_v)}
$$
$$
t^{n(f^-_b-f^+_b)}(1-t)^{f^-_r}(1-t^{-1})^{f^+_r}  t^{-(n-1)f^+_r}
\prod_{e\in F^+_r} t^{|P(f,e)|}
\prod_{e\in F_r}(n-|P(f,e)|)_t  =
$$
$$
t^{\d(K,n)}\sum_{f\in \F(G_K)}\mult_{t}(f)t^{\d(f)}t^{n(f^-_b-f^+_b)} 
(1-t)^{f^-_r}(1-t^{-1})^{f^+_r} 
\prod_{e\in F_r}(n-|P(f,e)|)_{t^{-\sign(s(e))}}=
$$
$$
t^{\d(K,n)}\sum_{f\in \F(G_K)}\mult_{t}(f)
t^{\d(f)} \prod_{v \in V_K} t^{-\sign(v)nf(e^b_v)}
\prod_{e \text{red};t(e)=v}\prod_{j=0}^{f(e)-1} 
(1-t^{-\sign(s(e))(n- j- \sum_{e'<_v e}f(e))}).  
$$


This finishes the proof.
\end{proof}

\appendix
\section{The zeta function of a graph and the Foata-Zeilberger formula}
\lbl{sec.FZ}

\subsection{The Foata-Zeilberger formula}
\lbl{sub.FZ}

In this section we translate key combinatorial results of 
Foata and Zeilberger \cite[Theorem 1.1]{FZ} in the language of our paper,
resulting in Theorem \ref{thm.FZ}.

Consider the {\em complete graph} $K_r$ with $r$ vertices equipped with a 
weight martix $B=(b_{ij})$ of size $r$ with independent commuting variables, 
and let $\R=\BZ[[b_{ij}]]$. Let $X=\{x_1,\dots,x_r\}$ denote an alphabet on 
$r$ letters and $X^\star$ denote the set of words on $X$.

Recall the notion of a {\em Lyndon word} $l \in X$, that is a word which is 
not a nontrivial power of another word, and is strictly smaller than any of 
its cyclic rearrangements. It follows by definition that

\begin{lemma}
\lbl{lem.lyndon}
There is a 1-1 correspondence between the set of nonperiodic cycles in $K_r$
and the set of Lyndon words in $X$.
\end{lemma}

Given a nonempty word $w=x_1 x_2 \dots x_m \in X$, Foata and Zeiberger define
a function $\b_{\circ}$ by 
$$
\b_{\circ}(w)=b_{x_1,x_2} b_{x_2, x_3} \dots b_{x_{m-1},x_m} b_{x_m,x_1}
$$
and $\b_{circ}(w)=1$ if $w$ is the empty word. 
Every word $w \in X$ has a {\em unique factorization} as
$$
w=l_1 l_2 \dots l_n
$$
where $l_i$ are Lyndon words in nonincreasing order 
$l_1 \geq l_2 \geq \dots \geq l_n$. Using this, Foata and Zeilberger define
a map:
$$
\b_{\dec}: X^\star \longto \R
$$
by
$\b_{\dec}(w)=\b_{\circ}(l_1) \b_{\circ}(l_2) \dots \b_{\circ}(l_n)$
where $(l_1,\dots,l_n)$ is the unique factorization of $w$. 
For example, if $X=\{1,2,3,4,5\}$ and $w=34512421231242$, then its 
factorization is given by $(l_1,l_2,l_3)=(345,1242,1231242)$ and
$\b_{\dec}(w)=b_{1,2}^3 b_{2,1}^2 b_{2,3} b_{2,4}^2 b_{3,1} b_{3,4} b_{4,2}^2
b_{4,5} b_{5,3}$.
 
Foata and Zeilberger define another map
$$
\b_{\vert}: X^\star \longto \R 
$$
as follows: if $w=x_1 x_2 \dots x_m$ is a word and $\ti{w}=\ti{x}_1 \ti{x}_2
\dots \ti{x}_m$ is the rearrangement of the letters of $w$ in nondecreasing 
order, then they define
$$
\b_{\vert}(w)=b_{\ti{x}_1,x_1} b_{\ti{x}_2,x_2} \dots b_{\ti{x}_m,x_m}. 
$$

In \cite[Theorem 1.1]{FZ} they show that
\begin{eqnarray}
\lbl{eq.FZ1}
\frac{1}{\det(I-B)} &=& \sum_{w \in X^\star} \b_{\dec}(w)\\ 
&=& \prod_{c \in \P(K_r)} \frac{1}{1-\b(c)}\\
\lbl{eq.FZ2}
&=& \sum_{w \in X^\star} \b_{\vert}(w) \in \R
\end{eqnarray}

Let us now translate \eqref{eq.FZ2}. Write a word $w$ and its rearrangement
$\ti w$ as an array
$$
\Ga(w)= \left[
\begin{matrix}
 \ti w \\
 w   
\end{matrix}
\right]
$$
A rearrangement $\ti w$ of a word $w$
is always of the form $\ti w=1^{n_1} 2^{n_2} \dots r^{n_r}$, and gives
rise to a function $f_w:\mathrm{Edges}(K_r)\longto \BN$ on the edges of $K_r$ 
defined by 
$f_w((i,j))$ is the number that the column vector $\left[
\begin{matrix}
 i \\
 j  
\end{matrix}
\right]
$ appears in $\Ga(w)$.
Since $\ti w$ is a rearrangement of $w$, it follows that $f_w$ is a flow. 
It follows from \ref{eq.FZ1} that this map $X^\star \longto \F(K_r)$
is onto, and it is easy to see that given an flow $\ga$ on $K_r$, the preimage
under this map consists of $\mult(\ga)$ words with the same $\b_{\vert}$
weight, equal to the weight of $\ga$. This together with
Equation \eqref{eq.FZ2} implies that
\begin{equation}
\lbl{eq.FZnew2}
\frac{1}{\det(I-B)}=
\sum_{f \in \F(K_r)} \b(f)\, \mult(f).
\end{equation}
This, together with a specialization of the variables imply Theorem 
\ref{thm.FZ}.

\section{A state sum for the Jones polynomial}
\lbl{sec.LW}

In this section we review the proof of Theorem \ref{thm.arcjones}.
The Jones polynomial $V$ of a link is determined by the skein theory:
$$
q^2 V\left(\psdraw{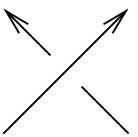}{.3in}\right)-q^{-2} V\left(\psdraw{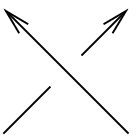}{.3in}\right)
=(q^{}-q^{-1}) V\left(\psdraw{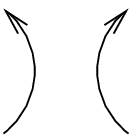}{.3in}\right)
$$
together with the initial condition $V(\unknot)(q)=q^{}+q^{-1}$. We will
be using a normalized version of the Jones polynomial defined by
$$ 
J(\K)(t)=V(\K)(t^{1/2})/V(\unknot)(t^{1/2}). 
$$

We review a state sum definition of the Jones polynomial
$V$ discussed by Turaev \cite{Tu} (see also \cite{J1}) and further
studied by Lin and Wang \cite{LW}. 
We recall the details of Turaev's general state sum construction, adapted
to our special case.

\begin{definition}
\lbl{def.state}
Fix a planar projection $\K$ of a knot. 
\newline
{\rm (a)}
Let $P_{\K}$ denote 
the planar digraph obtained from $\K$
by turning each crossing into a vertex. We call the edges of $P_{\K}$ 
{\em partarcs} of $\K$. 
\newline
{\rm (b)} 
A {\em state} $s$ of $\K$ is the assignment
of $0$ or $1$ to each partarc of $\K$, such that at each
crossing, the multiset of labels of the incoming edges equals
to the multiset of labels of outgoing edges. In other words, at each crossing
(positive or negative) a state looks like one of the following pictures,
$$
\psdraw{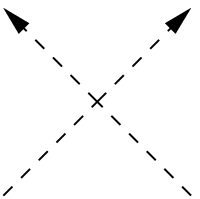}{0.5in} \hspace{0.5cm}
\psdraw{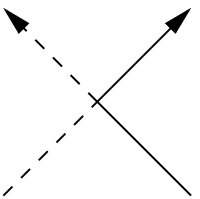}{0.5in} \hspace{0.5cm}
\psdraw{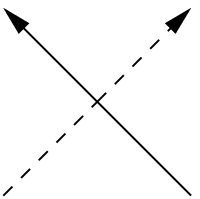}{0.5in} \hspace{0.5cm}
\psdraw{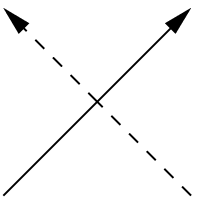}{0.5in} \hspace{0.5cm}
\psdraw{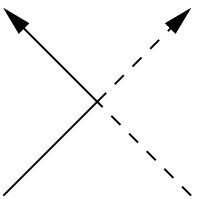}{0.5in} \hspace{0.5cm}
\psdraw{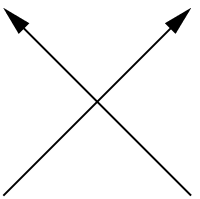}{0.5in} \hspace{0.5cm}
$$
where edges colored by $0$ or $1$ are depicted as dashed or solid 
respectively. 
\newline
{\rm (c)}
The {\em local weight} $\Pi_v(s)$ of a vertex $v$
of $P_{\K}$ of a state $s$ is given by
$$
\psdraw{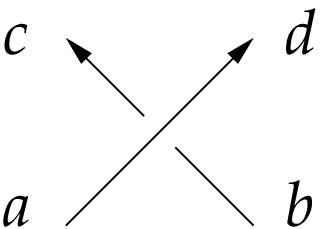}{0.5in}\longto (R^+)^{cd}_{ab} \hspace{1cm}
\psdraw{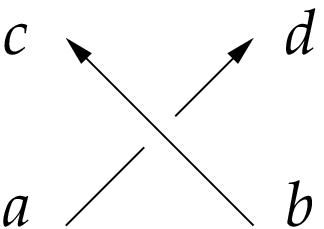}{0.5in}\longto (R^-)^{cd}_{ab}
$$
where $R^+$ and $R^-=(R^+)^{-1}$ is the $R$-matrix of the quantum group
$U_q(\mathfrak{sl}_2)$ given by:
\begin{align*}
(R^+)^{0,0}_{0,0}=(R^+)^{1,1}_{1,1}&=-q &
(R^+)^{1,0}_{0,1}=(R^+)^{0,1}_{1,0}&=1 &
(R^+)^{0,1}_{0,1}&=\bq-q  
\\
(R^-)^{0,0}_{0,0}=(R^-)^{1,1}_{1,1}&=-\bq &
(R^-)^{1,0}_{0,1}=(R^-)^{0,1}_{1,0}&=1 &
(R^-)^{1,0}_{1,0}&=q-\bq  
\end{align*}
where $\bq=q^{-1}$ and all other entries of the $R$ matrix are zero.
\newline
{\rm (d)}
The {\em weight} $\Pi(s)$ of a state $s$ is defined by 
$$
\Pi(s)=\prod_v \Pi_v(s).
$$  
Note that $(R^-)^{0,1}_{0,1}=(R^+)^{1,0}_{1,0}=0$. 
\newline
{\rm (e)}
A state $s$ is {\em admissible} iff $\Pi(s) \neq 0$. 
\end{definition}

There is an involution $s\to s^c$ of states of $\K$, obtained by interchanging
$0$ by $1$'s.

\begin{lemma}
\lbl{lem.arc}
{\rm (a)}
There is a 1-1 correspondence 
$$
\{ \, \text{\rm{states of}} \quad \K \} \longleftrightarrow
\{ \, \text{\rm{even subgraphs of}} \quad P_{\K} \}.
$$
{\rm (b)}
There is a 1-1 correspondence 
$$
\{ \, \text{\rm{admissible states of}} \quad \K \} \longleftrightarrow
\{ \, \text{\rm{collections of vertex-disjoint cycles of}} \quad G_{\K} \}.
$$
\end{lemma}

\begin{proof}
A state $s$ gives rise to an {\em even subgraph} of the 
$P_{\K}$ (whose edges are the ones colored by $1$ in $s$), also denoted by
$s$. Part (a) follows.

Since every vertex of $P_{\K}$ has outdegree $2$, it follows that the
involution of states corresponds to the operation of taking the
complement of an even subgraph in $P_{\K}$.
 
For part (b), 
observe that an admissible even subgraph $s$ of $P_{\K}$ gives rise 
to an even subgraph of the arc-graph $G_{\K}$ with each indegree at most one:
this follows since as mentioned above, $(R^+)^{0,1}_{0,1}=(R^-)^{1,0}_{1,0}=0$,
and so if we walk on $s$ along the orientation of $\K$, we never 
'jump down'; hence whenever we get to an arc of $\K$, we traverse it 
(along its orientation) until its end. Hence we can get to each arc at most 
once and $s$ corresponds 
to an even subgraph of $G_{\K}$ where each indegree is at most one.

Conversely, an even subgraph of $G_{\K}$ gives rise to a flow
on $P_{\K}$. This flow will be an admissible even subgraph of 
$P_{\K}$ if each indegree is at most one. The following figure illustrates the
excluded possibilities, where the value of the flow is shown on the
partarcs:
$$
\psdraw{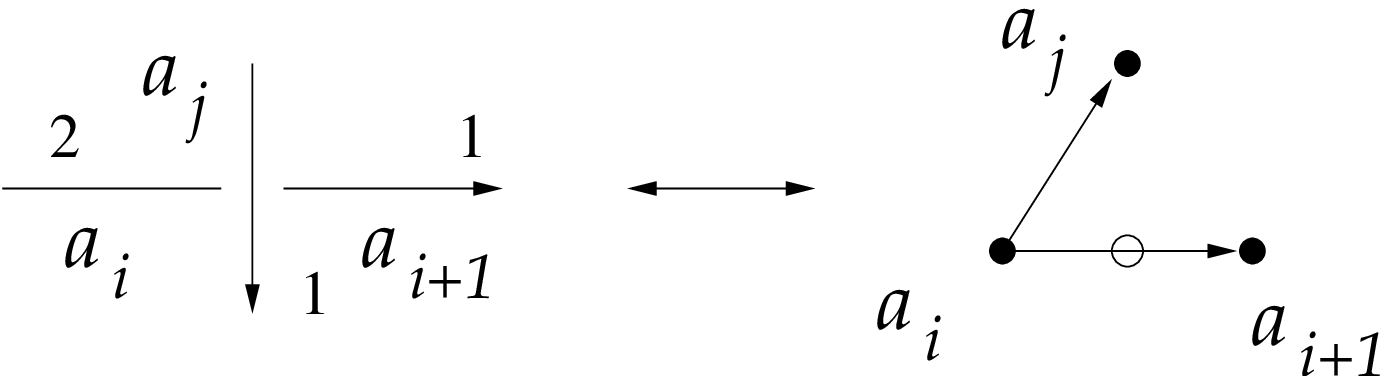}{3in}
$$
\end{proof}

\begin{definition}
\lbl{def.adm}
An even subgraph $G$ of $G_{\K}$ is {\em admissible} if each indegree 
is at most one. In other words, $G$ is a vertex-disjoint collection of 
directed cycles. Let $\Sev(G_{\K})$ denote the collection of 
admissible subgraphs of the arc-graph $G_{\K}$. 
\end{definition}

Next we define rotation and excess numbers of states.

\begin{definition}
\lbl{def.rottenp}
{\rm (a)} The {\em rotation number} $\rot(s)$ of a state $s$ 
is the number of counterclockwise circles colored by $1$
minus the number of clockwise circles colored by $1$, obtained from smoothening
of $s$, i.e., by the replacement:
$$
\psdraw{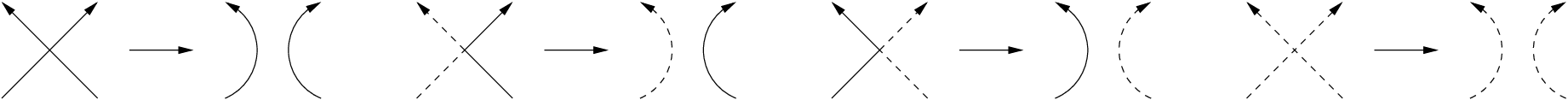}{4in}
$$
at all crossings of $s$.  
\newline
{\rm (b)} The {\em excess number}
$\exc(s)$ of a state $s$ is the sum of the signs of the crossings where
{\em all} four edges are colored by $1$ in $s$.
\end{definition}

With these preliminaries, the Jones polynomial is given by the state sum

$$
V(\K)(q)=(-q^2)^{-\w(\K)} \sum_{s \, \text{admissible}} 
q^{\rot(s^c)-\rot(s)} \Pi(s).
$$  

It was observed by Lin and Wang that the local weights of the $R$-matrix 
are proportional, up to a power of $q$ to the weights of a random walk
on $\K$. This is formalized in the following Lemma:

\begin{lemma}\cite[Lemma 2.3]{LW}
\lbl{lem.Rw}
For an admissible state $s$ of $\K$, we have:
$$
\Pi(s)=(-q)^{\w(\K)} q^{2 \exc(s)} \b(s)|_{t\to q^2}.
$$
\end{lemma}

\begin{proof}
First note that $\b(s)$ is well defined since by Lemma \ref{lem.arc} there
is a 1-1 correspondence between admissible states of $K$ and even admissible
subgraphs of $G_{\K}$, and each even subgraph is naturally a flow on $G_{\K}$.

Consider the following table of a state around a positive crossing:

\begin{center}
\begin{tabular}{|c|c|c|c|c|c|c|} \hline
 & $\psdraw{a1}{0.2in}$ & $\psdraw{a2}{0.2in}$ & $\psdraw{a3}{0.2in}$ & 
$\psdraw{a4}{0.2in}$ & $\psdraw{a5}{0.2in}$ & $\psdraw{a6}{0.2in}$ \\ \hline
$R$ & $-q$ & $\bar{q}-q$ & $1$ & $1$ & $0$ & $-q$ \\ \hline
$-\bar{q} R$ & $1$ & $1-\bar{q}^2$ & $-\bar{q}$ & $-\bar{q}$ & $0$ & $1$ 
\\ \hline
$\b$ & $1$ & $1-\bq^2$ & $\bq^2$ & $1$ & $0$ & $\bq^2$ \\ \hline
$q^{\exc}$ & $1$ & $1$ & $1$ & $1$ & $1$ & $q^2$ \\ \hline
$q^{\err}$ & $1$ & $1$ & $-q$ & $-\bq$ & $1$ & $1$ \\ \hline
\end{tabular}
\end{center}
and around a negative crossing: 

\begin{center}
\begin{tabular}{|c|c|c|c|c|c|c|} \hline
 & $\psdraw{a1}{0.2in}$ & $\psdraw{a2}{0.2in}$ & $\psdraw{a3}{0.2in}$ & 
$\psdraw{a4}{0.2in}$ & $\psdraw{a5}{0.2in}$ & $\psdraw{a6}{0.2in}$ \\ \hline
$R$ & $-\bq$ & $0$ & $1$ & $1$ & $q-\bq$ & $-\bq$ \\ \hline
$-q R$ & $1$ & $0$ & $-q$ & $-q$ & $1-q^2$ & $1$ \\ \hline
$\b$ & $1$ & $0$ & $1$ & $q^2$ & $1-q^2$ & $q^2$ \\ \hline
$q^{\exc}$ & $1$ & $1$ & $1$ & $1$ & $1$ & $\bq^2$ \\ \hline
$q^{\err}$ & $1$ & $1$ & $-q$ & $-\bq$ & $1$ & $1$ \\ \hline
\end{tabular}
\end{center}

Here, $\b(s)$ of a state $s$ equals to the weight of the $1$-part of $s$.

Inspection of these tables reveals that given a state $s$ and a crossing
of sign $\e=\pm 1$, we have $R=(-q)^{\e} q^{2 \exc+\err} \b$. Taking a product
over all vertices, we obtain that
$$
\Pi(s)=(-q)^{\w(\K)} q^{2 \exc(s)+\err(s)} \b(s)|_{t \to q^2}.
$$
It remains to show that $q^{\err(s)}=1$. $\err(s)$ is computed from a 
smoothening $\smooth(s)$ of $s$, which consists of a number of transversely
intersecting circles colored by $0$ or $1$. Any two transverse planar circles 
intersect on an even number of points, which can be paired up by paths 
on each circle. A case by case argument shows that
$\err(s)=1$. Some cases of the local contributions to $'\err'$ 
and their pairwise canceling is shown by:
$$
\psdraw{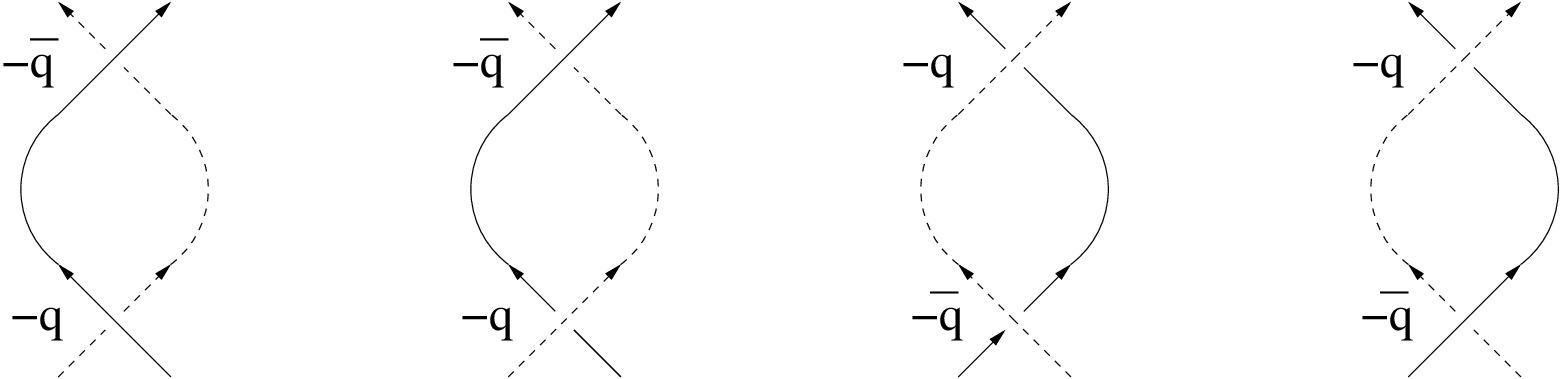}{5in}
$$
This concludes the proof of the lemma.
\end{proof}

The involution on the set of states of $\K$ has further consequences
discovered by Lin and Wang. Fix a partarc of $\K$ that borders the unbounded
region of the planar projection and mark it by $\star$. 
Let $\F'(\K-\star)$ denote the set of all admissible states
of $\K$ where $\star$ is colored by $0$.  

We will show first that
\begin{equation}
\lbl{eq.LW2}
J(\K)(t)= t^{\d(K)} \sum_{s \in \F'(K-\star)}  
t^{\d(s)} \b(s).
\end{equation}
We recall that
$ \d(K)=1/2(-\w(K)+\rot(K))$ and $\d(s)=\exc(s)-\rot(s)$.

Consider a long knot $K^{\text{long}}$ depicted as a box
and the two ways of closing it to obtain a knot $\K$ as follows:
$$
\psdraw{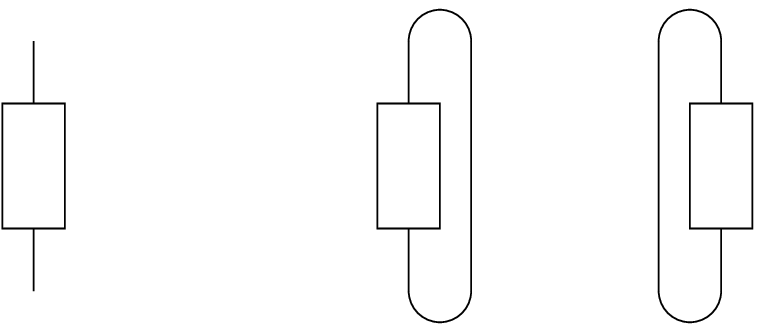}{1in}
$$
Let $a_i$ denote $V(K^{\text{long}})$ with boundary conditions $i$,
for $i=0,1$. Then, the two ways of closing $K^{\text{long}}$ give:
$$
V(K)(q)=q a_0 + q^{-1} a_1 = q^{-1} a_0 + q a_1
$$
from which follows that $a_0=a_1$ and thus $V(K)(q)=(q+q^{-1}) a_0=
V(\text{unknot}) a_0$. Thus, \eqref{eq.LW2} follows.

Next we introduce the rotation and excess numbers of a collection of
vertex disjoint cycles of $G_K$, using Lemma \ref{lem.arc}.

\subsection{Rotation and Excess numbers}
\lbl{sub.excess}

We observe that there is an integer function $\rot$  on the set of the edges
of $G_{\K}$ so that for each admissible state $s$ and its corresponding 
(see Lemma \ref{lem.arc})
 admissible subgraph $c$ of $G_{\K}$, $\rot(s)= \sum_{e\in c}\rot(e)$.

\begin{definition}
\lbl{def.jjj}
There is a Gauss map $d: \K \to S^1$ which together with the orientation of $\K$ and 
the counterclockwise orientation of $S^1$ induces a map
$$
H_1(P_{\K},\BZ) \longto H_1(S^1,\BZ)\cong \BZ.
$$
The above composition is defined to be the rotation number $\rot$. We can
think of the rotation number as an element of $H^1(P_{\K},\BZ)$
represented by a 1-cocyle, that is a map
$$
\rot: \Edges(P_{\K}) \longto \BZ.
$$
Consider now the arc-graph $G_{\K}$ of $\K$. 
There is a canonical map
$$
\Edges(G_{\K}) \longto 2^{\Edges(P_{\K})}
$$
defined as follows: if $(i,j)$ is an edge of $G_{\K}$, consider the $i$th
crossing of $P_{\K}$, and start walking on the part of the arc $a_j$ in
a direction of the orientation of $K$, until the end of the arc $a_j$.
This defines a collection of part-arcs that we associate to the edge $(i,j)$
of $G_{\K}$. Taking the sum of the rotation numbers of these part-arcs, 
defines a map 
$$
\rot: \Edges(G_K)\longto \BZ.
$$
\end{definition}

Next we show that $\exc'$ of next definition agrees
with $\exc$ of Definition \ref{def.rotten}.

\begin{definition}
\lbl{def.rottenp1}
Let $c$ be an admissible subgraph of $G_K$. We let $\exc'(c)$ equal
to $\exc(s)$, where $s$ is the corresponding admissible state
(see Definition \ref{def.rottenp} and Lemma \ref{lem.arc} for the 
correspondence).
\end{definition}

\begin{lemma}
\lbl{lem.exc1}
For every admissible subgraph $c$ of $G_K$, we have:
$$
\exc'(c)=\exc(c).
$$
\end{lemma}

\begin{proof}
$\exc'(c)$ is the sum of $\sign(v)$ where all 4 edges incident to a crossing
$v$ of $K$ belong to $c$:
$$
\psdraw{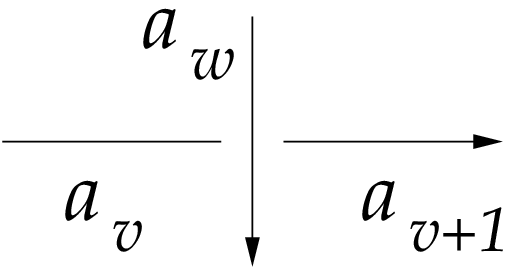}{1in}
$$ 
We will translate this in the language of the arc-graph $G_K$,
using Figure \ref{planarG2K}. A crossing $v$ as above determines a unique 
vertex of $G_K$ (corresponding to the arc $a_v$ ending at $v$) and a unique 
pair of edges $(e,e')$ of $G_K$: $e$ is the blue edge that starts at $v$,
and $e'$ is the unique edge of $c$ that ends in $w$ and signifies the
transition on the arc $a_w$. The result follows.
\end{proof}

\begin{proof}(of Theorem \ref{thm.arcjones})

Assume (after possibly changing the orientation of the knot, which does
not change the Jones polynomial) that we mark by $\star$ the last partarc of 
an arc of $\K$.
Lemma \ref{lem.arc}, \eqref{eq.LW2}, and subsection \ref{sub.excess}
 conclude the proof of Theorem \ref{thm.arcjones}.
\end{proof}

\section{A combinatorial counting of structures}
\lbl{sec.count}

In this section we consider structures on a set $[k]=\{1,\dots, k\}$,
and their combinatorial countings.

\begin{definition}
\lbl{def.str}
Let 
$k$ be a positive integer. A $k$-{\em structure} is a pair $S=(A,B)$
such that 
\begin{itemize}
\item
$A=(A_1,\dots,A_l)$, $B=(B_1,\dots,B_l)$ for some $l$, $A_i, B_i \subset [k]$,
$A_i \neq \emptyset$ for all $i$,
\item
$A$ is a partition of $\{1,\dots,k\}$ such that
for every $i < j$, $x\in A_i, y\in A_j$ we have $x < y$.
\item
$B_i\subset \cup_{j<i}A_j$. In particular $B_1=\emptyset$.
\item
$B$ is monotonic. That is, 
if $x\in B_j\cap A_i$ then for each $j\geq j'>i$, $x\in B_{j'}$,
\end{itemize}
\end{definition}


\begin{lemma}
\lbl{lem.numstr}
The number of $k$-structures $S$ such that $|A_i|=a_i$ and $|B_i|=b_i$
for $i=1,\dots,l$ 
is 
$$
\prod_{i=2}^l \binom{a_{i-1}+b_{i-1}}{b_i}.
$$
\end{lemma}

\begin{proof}
$B_i$ is an arbitrary subset of $A_{i-1}\cup B_{i-1}$ of $b_i$ elements. 
\end{proof}

\begin{definition}
\lbl{def.strdef}
Let $S$ be a $k$-structure, $v\in \{0,\dots, n-1\}^{\{1,\dots, k\}}$ and
$i\in A_x$ for some $x\leq l$.
\begin{itemize}
\item
We let $|S|=(a,b)$, where $a= (|A_1|,\dots, |A_l|)$ and 
$b= (|B_1|, \dots, |B_l|)$.
\item
We let $m(|S|,i)$ be the number of $j\in A_x\cup B_x$ such that $j<i$. 
Note that
$m(|S|,i)$ equals $b_i$ plus the number of elements of $A_x$ that are smaller 
than $i$ and hence it depends only on $|S|$.
\item
We denote by $\defect_1(S,v,i)$ the number of $j\in A_x\cup B_x$ such that 
$j<i$ and $v_j<v_i$,
\item
we denote by $\defect_2(S,v,i)$ the number of $j\in B_{d(i)+1}$ such that 
$v_j<v_i$; $d(i)$ is 
the minimum index so that $d(i)\geq i$ and $i\notin B_{d(i)+1}$.
\end{itemize}
\end{definition}

\begin{definition}
\lbl{def.strvec}
Let $S$ be a $k$-structure. 
We let $V(S,n)=\{v\in \{0,\dots, n-1\}^{\{1,\dots, k\}}; $
if $\{i,j\}\subset A_m\cup B_m$ for some $m$ then $v_i\neq v_j\}$. 
\end{definition}

The following Theorem follows by comparing the definitions.

\begin{theorem}
\lbl{thm.sorstr}
Let $f$ be a flow on arc-graph $G_K$. Recall that $f_r(v)=\sum f(e)$ over 
all red edges
of $G_K$ terminating in vertex $v$ of $G_K$, and $f_b(v)$ is defined 
analogously for the blue edges. 
We consider set $F_r$ linearly ordered, first by the terminal vertices, 
and then by ordering $\prec$
which induces a linear ordering on the set $\cup F(e)$, over red edges $e$ 
entering the same vertex
(see Definition \ref{def.orderf}).

There is a natural bijection between $\AC_n(f)$ and the set of all pairs 
$(S,v)$ where $S$ 
is an $|F_r|$-structure, $|S|=((f_r(1),\dots, f_r(r))(f_b(1),\dots, f_b(r)))$ 
and $v\in V(S,n)$.
\end{theorem}

\begin{theorem}
\lbl{thm.str10}
$$
\sum_{S: |S|=(a,b)} 
\sum_{v\in V(S,n)}\prod_{i=1}^k t^{v_i-\defect_1(S,v,i)-\defect_2(S,v,i)}
= \prod_{i=1}^k(n-m(i))_t
\prod_{i=1}^{l-1} \binom{a_i+b_i}{b_{i+1}}_{t^{-1}}.
$$
\end{theorem}

In the proof we will use the following proposition.

\begin{proposition}
\lbl{prop.str2}
Let $S$ be a $k$-structure. Then
$$
\sum_{v\in V(S,n)}\prod_{i=1}^k t^{v_i-\defect_1(S,v,i)}
= \prod_{i=1}^k (n-m(|S|,i))_{t}.
$$
\end{proposition}

\begin{proof}
Use induction on $k$. The inductive step follows from the following claim:

{\em Claim:}
Let $m(k)< n$ and fix different numbers $v_1,\dots, v_{m(k)}$ between $0$ 
and $n-1$.
Then 
\begin{itemize}
\item
$$
\sum_{v_k: v_k\neq v_i, i\leq m(k)} t^{v_k-\defect_1(v,k)}= A-B+C,
$$
where
$A=\sum_{v_k: v_k\neq v_i, i\leq m(k)} t^{v_k}$, $B=\sum_{i=1}^{m(k)} 
t^{n-i}$, and
$C=\sum_{i=1}^{m(k)} t^{v_i}$.
\item
$A+C= \sum_{0\leq z\leq n-1}t^z$ and $A-B+C= \frac{1-t^{n-m(k)}}{1-t}$.
\end{itemize}

Note that the second part is simply true. 

Let $v'_1 < \dots < v'_{m(k)}$ be a reordering of $v_1,\dots, v_{m(k)}$.
We may write $v'_1=n-i_1, \dots,v'_{m(k)}=n-i_{m(k)}, 
1\leq i_{m(k)}<\dots <i_1$. 
The LHS becomes
$$
A -t^{n-i_1+1}-\dots -t^{n-i_2-1}-t^{n-i_2+1}-\dots 
-t^{n-i_{m(k)}-1}-t^{n-i_{m(k)}+1}-\dots -t^{n-1}
$$
$$
+ t^{n-i_1}+\dots +t^{n-i_2-2}+t^{n-i_2-1}+\dots +t^{n-m(k)-1}.
$$
This equals to the RHS of the equality we wanted to show. The Proposition 
simply follows
from the Claim.
\end{proof}

\begin{proof}(of Theorem \ref{thm.str10})

We let $a'_i= \sum_{j\leq i}a_i$.
$$
 \sum_{S: |S|=(a,b)} \sum_{v\in V(S,n)}
\prod_{i=1}^k t^{v_i-\defect_1(S,v,i)-\defect_2(S,v,i)}=
$$
$$
\sum_{B_2\subset A_1}\sum_{v_1,\dots, v_{a_1}}
\prod_{i=1}^{a_1} t^{v_i-\defect_1(v,i)-\defect_2(B_2,v,i)}\times
\sum_{B_3, \dots, B_l}\sum_{v_{a_1+1},\dots, v_k}
\prod_{i=a_1+1}^k t^{v_i-\defect_1(v,i)-\defect_2(B_3,\dots, B_l,v,i)}=
$$
$$
\sum_{v_1,\dots, v_{a_1}}\prod_{i=1}^{a_1} t^{v_i-\defect_1(v,i)}
\sum_{B_2\subset A_1}\prod_{i=1}^{a_1} t^{-\defect_2(B_2,v,i)}\times
\sum_{B_3, \dots, B_l}\sum_{v_{a_1+1},\dots, v_k}
\prod_{i=a_1+1}^{k} t^{v_i-\defect_1(v,i)-\defect_2(B_3,\dots, B_l,v,i)}=
$$
$$
\sum_{v_1,\dots, v_{a_1}}\prod_{i=1}^{a_1} t^{v_i-\defect_1(v,i)}
\sum_{B_2\subset A_1}\prod_{i=1}^{a_1} t^{-\defect_2(B_2,v,i)}\times
\dots \times
$$
$$
\sum_{v_{a'_{l-2}+1},\dots, v_{a'_{l-1}}}\prod_{i=a'_{l-2}+1}^{a'_{l-1}} 
t^{v_i-\defect_1(v,i)}
\sum_{B_l}\prod_{i=a'_{l-2}+1}^{a'_{l-1}}t^{-\defect_2(B_l,v,i)}\times
\sum_{v_{a'_{l-1}+1},\dots, v_k}
\prod_{i=a'_{l-1}+1}^k t^{v_i-\defect_1(v,i)}.
$$

The last sum may be expressed using Proposition \ref{prop.str2}, and we get

$$
\sum_{v_1,\dots, v_{a_1}}\prod_{i=1}^{a_1} t^{v_i-\defect_1(v,i)}
\sum_{B_2\subset A_1}\prod_{i=1}^{a_1} t^{-\defect_2(B_2,v,i)}\times
\dots \times
$$
$$
\sum_{v_{a'_{l-2}+1},\dots, v_{a'_{l-1}}}\prod_{i=a'_{l-2}+1}^{a'_{l-1}} 
t^{v_i-\defect_1(v,i)}
\sum_{B_l}\prod_{i=a'_{l-2}+1}^{a'_{l-1}}t^{-\defect_2(B_l,v,i)}\times
\prod_{i=a'_{l-1}+1}^k(n-m(i))_t=
$$
$$
\prod_{i=a'_{l-1}+1}^k(n-m(i))_t
\sum_{v_1,\dots, v_{a_1}}\prod_{i=1}^{a_1} t^{v_i-\defect_1(v,i)}
\sum_{B_2\subset A_1}\prod_{i=1}^{a_1} t^{-\defect_2(B_2,v,i)}\times
\dots \times
$$
$$
\prod_{i=a'_{l-2}+1}^{a'_{l-1}} (n-m(i))_t 
\binom{a_{l-1}+b_{l-1}}{b_l}_{t^{-1}}=
$$
$$
 \prod_{i=1}^k(n-m(i))_t\prod_{i=1}^{l-1} \binom{a_i+b_i}{b_{i+1}}_{t^{-1}}.
$$
\end{proof}

\ifx\undefined\bysame
        \newcommand{\bysame}{\leavevmode\hbox
to3em{\hrulefill}\,}
\fi

\end{document}